\documentclass[11pt]{article}
\usepackage{graphicx}
\usepackage{latexsym, amsmath, amsfonts, amssymb,amsthm}
\usepackage{mathptmx}
\usepackage[T1]{fontenc}
\usepackage[latin9]{inputenc}
\usepackage{color}
\usepackage{textcomp}
\usepackage{longtable}
\usepackage[english]{babel}
\providecommand{\U}[1]{\protect\rule{.1in}{.1in}}

\numberwithin{equation}{section} \numberwithin{figure}{section}
\theoremstyle{plain} \theoremstyle{plain}
\newtheorem{theorem}{Theorem}
\theoremstyle{definition}
\newtheorem{definition}[theorem]{Definition}
\theoremstyle{plain}
\newtheorem{lemma}[theorem]{Lemma}
\theoremstyle{plain}
\newtheorem{proposition}[theorem]{Proposition}
\theoremstyle{definition}
\newtheorem{example}[theorem]{Example}
\theoremstyle{plain}
\newtheorem{corollary}[theorem]{Corollary}
\theoremstyle{remark}

\newcommand{\gph}{\mbox{\rm Graph}}

\newcommand{\sur}{\mbox{\rm sur}\,}
\newcommand{\Sur}{\mbox{\rm Sur}\,}

\begin{document}

\date{\today}


\begin{center}
{\Large \textbf{Continuity of set-valued maps
\\ \smallskip revisited in the light of tame geometry}}

\vspace{0.8cm}

{\large \textsc{Aris Daniilidis \& C. H. Jeffrey Pang}}
\end{center}

\bigskip

\noindent\textbf{Abstract} Continuity of set-valued maps is hereby
revisited: after recalling some basic concepts of varia\-tional
analysis and a short description of the State-of-the-Art, we obtain
as by-product two Sard type results concerning local minima of
scalar and vector valued functions. Our main result though, is
inscribed in the framework of tame geometry, stating that a
closed-valued semialgebraic set-valued map is almost everywhere
continuous (in both topological and measure-theoretic sense). The
result --depending on stratification techniques-- holds true in a
more general setting of o-minimal (or tame) set-valued maps. Some
applications are briefly discussed at the end.

\bigskip


\bigskip


\bigskip

\tableofcontents

\bigskip

\section{Introduction}

We say that $S$ is a \emph{set-valued map} (we also use the term
\emph{multivalued function} or simply \emph{multifunction}) from $X$
to $Y$, denoted by $S:X\rightrightarrows Y$, if for every $x\in X$,
$S\left( x\right)  $ is a subset of $Y$. All single-valued maps in
classical analysis can be seen as set-valued maps, while many
problems in applied mathematics are set-valued in nature. For
instance, pro\-blems of stability (parametric optimization) and
controllability are often best treated with set-valued maps, while
gradients of (differentiable) functions, tangents and normals of
sets (with a structure of differentiable manifold) have natural
set-valued generalizations in the nonsmooth case, by means of
variational analysis techniques. The inclusion $y\in S(x)$ is the
heart of modern variational analysis. We refer the reader to
\cite{AF90,RW98} for more details.\smallskip

Continuity properties of set-valued maps are crucial in many
applications. A typical set-valued map arising from some
construction or variational problem will not be continuous.
Nonetheless, one often expects a kind of semicontinuity (inner or
outer) to hold. (We refer to Section 2 for relevant
definitions.)\smallskip

A standard application of a Baire argument entails that
closed-valued set-valued maps are generically continuous, provided
they are either inner or outer semicontinuous. Recalling briefly
these results, as well as other concepts of continuity for
set-valued maps, we illustrate their sharpness by means of
appropriate examples. We also mention an interesting consequence of
these results by establishing a Sard-type result for the image of
local minima.

Moving forward, we limit ourselves to semialgebraic maps
\cite{BR90,Coste99} or more generally, to maps whose graph is a
definable set in some o-minimal structure \cite{vdDries98, Coste02}.
This setting aims at eliminating most pathologies that pervade
analysis which, aside from their indisputable theoretical interest,
do not appear in most practical applications. The definition of a
definable set might appear reluctant at the first sight (in
particular for researchers in applied mathematics), but it
determines a large class of objects (sets, functions, maps)
encompassing for instance the well-known class of semialgebraic sets
\cite{BR90,Coste99}, that is, the class of Boolean combinations of
subsets of $\mathbb{R}^{n}$ defined by finite polynomials and
inequalities. All these classes enjoy an important stability
property ---in the case of semialgebraic sets this is expressed by
the Tarski-Seidenberg (or quantifier elimination) principle--- and
share the important property of stratification: every definable set
(so in particular, every semialgebraic set) can be written as a
disjoint union of smooth manifolds which fit each other in a regular
way (see Theorem~\ref{theorem:stratification} for a precise
statement). This tame behaviour has been already exploited in
various ways in variational analysis, see for instance \cite{AB2009}
(convergence of proximal algorithm), \cite{BDL04} ({\L}ojasiewicz
gradient inequality), \cite{BDL2008} (semismoothness),
\cite{Ioffe08} (Sard-Smale type result for critical values) or
\cite{Ioffe2009} for a recent survey of what is nowadays called
\emph{tame optimization.}\smallskip

The main result of this work is to establish that every
semialgebraic (more generally, definable) closed-valued set-valued
map is generically continuous. Let us point out that in this
semialgebraic context, genericity implies that possible failures can
only arise in a set of lower dimension, and thus is equivalent to
the measure-theoretical notion of \emph{almost-everywhere} (see
Proposition~\ref{Prop:generic} for a precise statement). The proof
uses properties of stratification, some technical lemmas of
variational analysis and a recent result of
Ioffe~\cite{Ioffe08}.\smallskip

The paper is organized as follows. In Section~\ref{sec-2} we recall
basic notions of variational analysis and revisit results on the
continuity of set-valued maps. As by-product of our development we
obtain, in Section~\ref{sec-3} two Sard-type results: the first one
concerns minimum values of (scalar) functions, while the second one
concerns Pareto minimum values of set-valued maps. We also grind our
tools by adapting the Mordukhovich criterion to set-valued maps with
domain a smooth submanifold ${\mathcal{X}}$ of $\mathbb{R}^{n}$. In
Section~\ref{sec-main} we move into the semialgebraic case. Adapting
a recent result of Ioffe \cite[Theorem~7]{Ioffe08} to our needs, we
prove an intermediate result concerning generic strict continuity of
set-valued maps with a closed semialgebraic graph. Then, relating
the failure of continuity of the mapping with the failure of its
trace on a stratum of its graph, and using two technical lemmas we
establish our main result. Section~\ref{sec-appl} contains some
applications of the main result.

\bigskip

\textbf{Notation.} Denote $\mathbb{B}^{n}\left(x,\delta\right)$ to
be the closed ball of center $x$ and radius $r$ in $\mathbb{R}^{n}$,
and $\mathbb{S}^{n-1}\left(x,r\right)$ to be the sphere of center
$x$ and radius $r$ in $\mathbb{R}^{n}$. When there is no confusion
of the dimensions of $\mathbb{B}^{n}\left( x,r\right)$ and
$\mathbb{S}^{n-1}\left( x,r\right)$, we omit the superscript. The
unit ball $\mathbb{B}\left( \mathbf{0},1\right)$ is denoted by
$\mathbb{B}$. We denote by $\mathbf{0}_n$ the neutral element of
$\mathbb{R}^n$. As before, if there is no confusion on the dimension
we shall omit the subscript. Given a subset $A$ of $\mathbb{R}^{n}$
we denote by $\mathrm{cl\,}(A)$, $\mathrm{int}(A)$ and $\partial A$
respectively, its topological closure, interior and boundary. For
$A_{1} ,A_{2}\subset\mathbb{R}^{n}$ and $r\in\mathbb{R}$ we set
\[
A_{1}+rA_{2}:=\{a_{1}+ra_{2}:a_{1}\in A_{1},a_{2}\in A_{2}\}.
\]
We recall that the Hausdorff distance $\mathcal{D}(A_{1},A_{2})$
between two bounded subsets $A_{1},A_{2}$ of $\mathbb{R}^{n}$ is
defined as the infimum of all $\delta>0$ such that both inclusions
$A_{1}\subset A_{2}+\delta \,\mathbb{B}$ and $A_{2}\subset
A_{1}+\delta\,\mathbb{B}$ hold (see \cite[Section~9C]{RW98} for
example). Finally, we denote by $$\mathrm{Graph }(S)=\{(x,y)\in
X\times Y: \, y\in S(x)\}\,,$$ the graph of the set-valued map
$S:X\rightrightarrows Y$.

\section{Basic notions in set-valued analysis}

\label{sec-2}

In this section we recall the definitions of continuity (outer, inner, strict)
for set-valued maps, and other related notions from variational analysis. We
refer to \cite{AF90,RW98} for more details.

\subsection{Continuity concepts for set-valued maps}

We start this section by recalling the definitions of continuity for
set-valued maps.

\medskip

\textbf{(Kuratowski limits of sequences)} We first recall basic notions about
(Kuratowski) limits of sets. Given a sequence $\left\{  C_{\nu}\right\}
_{\nu\in\mathbb{N}}$ of subsets of $\mathbb{R}^{n}$ we define:

\begin{itemize}
\item the \emph{outer limit} $\limsup_{\nu\rightarrow\infty}C_{\nu}$, as the
set of all accumulation points of sequences
$\{x_{\nu}\}_{\nu\in\mathbb{N} }\subset\mathbb{R}^{n}$ with
$x_{\nu}\in C_{\nu}$ for all $\nu\in\mathbb{N}$. In other words,
$x\in\limsup_{\nu\rightarrow\infty}C_{\nu}$ if and only if for every
$\varepsilon>0$ and $N\geq1$ there exists $\nu\geq N$ with
$C_{\nu}\cap B(x,\varepsilon)\neq\emptyset$ ;

\item the \emph{inner limit} $\liminf_{\nu\rightarrow\infty}C_{\nu}$, as the
set of all limits of sequences $\{x_{\nu}\}_{\nu\in\mathbb{N}}\subset
\mathbb{R}^{n}$ with $x_{\nu}\in C_{\nu}$ for all $\nu\in\mathbb{N}$. In other
words, $x\in\liminf_{\nu\rightarrow\infty}C_{\nu}$ if and only if for every
$\varepsilon>0$ there exists $N\in\mathbb{N}$ such that for all $\nu\geq N$ we
have $C_{\nu}\cap B(x,\varepsilon)\neq\emptyset.$
\end{itemize}

Furthermore, we say that the\emph{ limit }of the sequence $\left\{
C_{\nu }\right\}  _{\nu\in\mathbb{N}}$ exists if the outer and inner
limit sets are equal. In this case we write:
\[
\lim_{\nu\rightarrow\infty}C_{\nu}:=\limsup_{\nu\rightarrow\infty}C_{\nu
}=\liminf_{\nu\rightarrow\infty}C_{\nu}.
\]

\smallskip

\textbf{(Outer/inner continuity of a set-valued map) }Given a
set-valued map $S:\mathbb{R}^{n}\rightrightarrows\mathbb{R}^{m}$, we
define the outer (respectively, inner) limit of $S$ at
$\bar{x}\in\mathbb{R}^{n}$ as the union of all upper limits
$\limsup_{\nu\rightarrow\infty}S\left(  x_{\nu}\right)  $
(respectively, intersection of all lower limits
$\liminf_{\nu\rightarrow \infty}S\left(  x_{\nu}\right)  $) over all
sequences $\left\{  x_{\nu }\right\}  _{\nu\in\mathbb{N}}$
converging to $\bar{x}$. In other words:
\[
\limsup_{x\rightarrow\bar{x}}\,S\left(  x\right) \,
:=\bigcup_{x_{\nu}
\rightarrow\bar{x}}\limsup_{\nu\rightarrow\infty}\,S\left(
x_{\nu}\right)
\qquad\text{and}\qquad\liminf_{x\rightarrow\bar{x}}\,S\left(
x\right)
:=\bigcap_{x_{\nu}\rightarrow\bar{x}}\liminf_{\nu\rightarrow\infty}\,S\left(
x_{\nu}\right)  .
\]

We are now ready to recall the following definition.

\begin{definition}
\label{def:osc-isc-map}\cite[Definition~5.4]{RW98} A set-valued map
$S:\mathbb{R}^{n}\rightrightarrows\mathbb{R}^{m}$ is called
\emph{outer semicontinuous} at $\bar{x}$ if
\[
\limsup_{x\rightarrow\bar{x}}\,S\left(  x\right)  \subset S\left(
\bar {x}\right)\,,
\]
or equivalently, $\limsup_{x\rightarrow\bar{x}}\,S\left(  x\right)
=S\left( \bar{x}\right)  $, and \emph{inner semicontinuous} at
$\bar{x}$ if
\[
\liminf_{x\rightarrow\bar{x}}\,S\left(  x\right)  \supset S\left(
\bar {x}\right)  ,
\]
or equivalently when $S$ is closed-valued,
$\liminf_{x\rightarrow\bar{x} }S\left(  x\right)  =S\left(
\bar{x}\right)  $. It is called \emph{continuous} at $\bar{x}$ if
both conditions hold, i.e., if $S\left( x\right)  \rightarrow
S\left(  \bar{x}\right)  $ as $x\rightarrow\bar{x}$.

If these terms are invoked relative to $X$, a subset of $\mathbb{R}^{n}$
containing $\bar{x}$, then the properties hold in restriction to convergence
$x\rightarrow\bar{x}$ with $x\in X$ (in which case the sequences $x_{\nu
}\rightarrow\bar{x}$ in the limit formulations are required to lie in~$X$).
\end{definition}

Notice that every outer semicontinuous set-valued map has closed
values. In particular, it is well known that

\begin{itemize}
\item $S$ is outer semicontinuous if and only if $S$ has a closed graph.
\end{itemize}

When $S$ is a single-valued function, both outer and inner
semicontinuity reduce to the standard notion of continuity. The
standard example of the mapping
\begin{equation}
S\left(  x\right)  :=
\begin{cases}
0 & \mbox{if }x\mbox{ is rational}\\
1 & \mbox{if }x\mbox{ is irrational}
\end{cases}
\label{eq:rat-irrat}
\end{equation}
shows that it is possible for a set-valued map to be nowhere outer
and nowhere inner semicontinuous. Nonetheless, the following
genericity result holds. (We recall that a set is \emph{nowhere
dense} if its closure has empty interior, and \emph{meager} if it is
the union of countably many sets that are nowhere dense in $X$.) The
following result appears in \cite[Theorem~5.55]{RW98} and
\cite[Theorem~1.4.13]{AF90} and is attributed to
\cite{K32,Cho47,Shi81}. The domain of $S$ below can be taken to be a
complete metric space, while the range can be taken to be a complete
separable metric space, but we shall only state the result in the
finite dimensional case.

\begin{theorem}
\label{theorem:[RW98,5.55]} Let $X\subset \mathbb{R}^{n}$ and
$S:\mathbb{R}^{n}\rightrightarrows\mathbb{R}^{m}$ be a closed-valued
set-valued map. Assume $S$ is either outer semicontinuous or inner
semicontinuous relative to $X$. Then the set of points $x\in X$
where $S$ fails to be continuous relative to $X$ is meager in $X$.
\end{theorem}

The following example shows the sharpness of the result, if we move
to incomplete spaces.

\smallskip

\begin{example}
Let $c_{00}(\mathbb{N})$ denote the vector space of all real
sequences $x=\{x_{n}\}_{n\in\mathbb{N}}$ with finite support
$\mathrm{supp} (x):=\{i\in\mathbb{N}:x_{i}\neq0\}$. Then the
operator $S_{1}(x)=$ $\mathrm{supp}(x)$ is everywhere inner
semicontinuous and nowhere outer semicontinuous, while the operator
$S_{2}(x)=\mathbb{Z}\setminus S_{1}(x)$ is everywhere outer
semicontinuous and nowhere inner semicontinuous.$\hfill\Box$
\end{example}

\medskip

\textbf{(Strict continuity of set-valued maps) }A stronger concept
of continuity for set-valued maps is that of \emph{strict
continuity} \cite[Definition~9.28]{RW98}, which is equivalent to
Lipschitz continuity when the map is single-valued. For set-valued
maps $S:\mathbb{R} ^{n}\rightrightarrows\mathbb{R}^{m}$ with bounded
values, strict continuity is quantified by the Hausdorff distance.
Namely, a set-valued map $S$ is strictly continuous at $\bar{x}$
(relative to $X$) if the quantity
\[
\mathrm{lip}_{X}S(\bar{x}):=\underset{{\scriptsize
\begin{array}
[c]{c}
x,x^{\prime}\rightarrow\bar{x}\\
x\neq x^{\prime}
\end{array}
}}{\lim\sup}\,\frac{\mathcal{D}\,(S(x),S(x^{\prime}))}{|x-x^{\prime}|}
\]
is bounded. In the general case (that is, when $S$ maps to unbounded
sets), we say that $S$ is strictly continuous, whenever the
truncated map $S_{r} :\mathbb{R}^{n}\rightrightarrows\mathbb{R}^{m}$
defined by
\[
S_{r}\left(  x\right)  :=S\left(  x\right)  \cap r\,\mathbb{B},
\]
is Lipschitz continuous for every $r>0$.

\subsection{Normal cones, coderivatives and the Aubin property}

Before we consider other concepts of continuity of set-valued maps
we need to recall some basic concepts from variational analysis. We
first recall the definition of the Hadamard and limiting normal
cones.

\begin{definition}
\label{def:Normal-cones} (Normal cones) \cite[Definition 6.3]{RW98} For a
closed set $D\subset\mathbb{R}^{n}$ and a point $\bar{z}\in D$, we recall that
the \emph{Hadamard normal cone} $\hat{N}_{D}\left(  \bar{z}\right)  $ and the
\emph{limiting normal cone} $N_{D}\left(  \bar{z}\right)  $ are defined by
\begin{align}
\hat{N}_{D}\left(  \bar{z}\right)   &  :=\left\{  v\mid\left\langle
v,z-\bar{z}\right\rangle \leq o\left(  \left\vert z-\bar{z}\right\vert
\right)  \mbox{ for }z\in D\right\}, \nonumber\\
N_{D}\left(  \bar{z}\right)   &  :=\{v\mid\exists\left\{  z_{i},v_{i}\right\}
_{i=1}^{\infty}\subset\mbox{\rm Graph}(\hat{N}_{D}),\;\nu_{i}\rightarrow
v\mbox{ and }z_{i}\rightarrow\bar{z}\}\nonumber\\
&  =\limsup_{z\rightarrow\bar{z},z\in D}\,\hat{N}_{D}\left(
z\right). \nonumber
\end{align}
When $D$ is a smooth manifold, both notions of normal cone coincide
and define the same subspace of $\mathbb{R}^{n}$. A dual concept to
the normal cone is the \emph{tangent cone} {$T_{D}\left(
\bar{z}\right)$. While tangent cones can be defined for nonsmooth
sets, our use here shall be restricted only to tangent cones of
manifolds, that is, tangent spaces in the sense of differential
geometry, in which case $T_{D}\left( \bar{z}\right)  =\left(
N_{D}\left(  \bar{z}\right)  \right)^{\perp}$.}
\end{definition}

As is well-known, the generalization of the adjoint of a linear operator for
set-valued maps is derived from the normal cones of its graph.

\begin{definition}
(Coderivatives) \cite[Definition~8.33]{RW98} For $F:\mathbb{R}^{n}
\rightrightarrows\mathbb{R}^{m}$ and $\left(  \bar{x},\bar{y}\right)
\in\mbox{\rm Graph}\left(  F\right)  $, the \emph{limiting
coderivative} $D^{\ast}F\left(  \bar{x}\mid\bar{y}\right)
:\mathbb{R}^{m}\rightrightarrows \mathbb{R}^{n}$ is defined by
\[
D^{\ast}F\left(  \bar{x}\mid\bar{y}\right)  \left(  y^{\ast}\right)  =\left\{
x^{\ast}\mid\left(  x^{\ast},-y^{\ast}\right)  \in
N_{{\scriptsize \mbox{Graph}\left(  F\right)  }}\left(  \bar{x},\bar
{y}\right)  \right\}  .
\]
It is clear from the definitions that the coderivative is a positively
homogeneous map, which can be measured with the outer norm below.

\begin{definition}
\cite[Section~9D]{RW98} The \emph{outer norm} $\left\vert
\cdot\right\vert ^{+}$ of a positively homogeneous map
$H:\mathbb{R}^{n}\rightrightarrows \mathbb{R}^{m}$ is defined by
\[
\left\vert H\right\vert ^{+}:=\sup_{w\in\mathbb{B}^{n}\left(
\mathbf{0} ,1\right)  }\sup_{z\in H\left(  w\right)  }\left\vert
z\right\vert =\sup\left\{  \frac{\left\vert z\right\vert
}{\left\vert w\right\vert } \mid\left(  w,z\right)  \in\mbox{\rm
Graph}\left(  H\right)  \right\}  .
\]
\end{definition}
\end{definition}

\smallskip

\textbf{(Aubin property and Mordukhovich criterion)} We now recall the Aubin
Property and the graphical modulus, which are important to study local
Lipschitz continuity properties of a set-valued map.

\begin{definition}
(Aubin property and graphical modulus) \cite[Definition 9.36]{RW98}
A map $S:\mathbb{R}^{n}\rightrightarrows\mathbb{R}^{m}$ has the
\emph{Aubin property relative to $X$ at $\bar{x}$ for $\bar{u}$},
where $\bar{x}\in X\subset\mathbb{R}^n$ and $\bar {u}\in S\left(
\bar{x}\right) $, if $\mbox{\rm Graph}\left(  S\right)  $ is locally
closed at $\left( \bar{x},\bar{u}\right)  $ and there are
neighborhoods $V$ of $\bar{x}$ and $W$ of $\bar{u}$, and a constant
$\kappa \in\mathbb{R}_{+}$ such that
\[
S\left(  x^{\prime}\right)  \cap W\subset S\left(  x\right)  +\kappa\left\vert
x^{\prime}-x\right\vert \mathbb{B}\mbox{ for all }x,x^{\prime}\in X\cap V.
\]
This condition with $V$ in place of $X\cap V$ is simply the
\emph{Aubin property at $\bar{x}$ for $\bar{u}$}. The
\emph{graphical modulus of $S$ relative to $X$ at $\bar{x}$ for
$\bar{u}$} is then
\begin{align*}
\mbox{\rm lip}_{X}S\left(  \bar{x}\mid\bar{u}\right)   &  := \inf\{\kappa
\mid\exists \mbox{ neighborhoods }V\mbox{ of }\bar{x}\mbox{ and }W \mbox{ of }\bar{u} \mbox{ s.t. }\\
& \qquad\qquad   S\left(  x^{\prime}\right)  \cap W\subset S\left(  x\right)  +\kappa
\left\vert x^{\prime}-x\right\vert \mathbb{B}\mbox{ for all }x,x^{\prime}\in
X\cap V\}.
\end{align*}
{In the case where $X=\mathbb{R}^n$, the subscript $X$ is omitted.}
\end{definition}

The following result (known as Mordukhovich criterion
\cite[Theorem~9.40]{RW98}) characterizes the Aubin property by means
of the corresponding coderivative. (For a primal characterization
using the graphical derivative see \cite[Theorem~1.2]{DQZ06}.)

\begin{proposition}
[Mordukhovich criterion]\label{pro:Boris-criterion}Let $S:\mathbb{R}
^{n}\rightrightarrows\mathbb{R}^{m}$ be a set-valued map whose graph
$\mbox{\rm Graph}\left(  S\right)  $ is locally closed at $\left(
\bar {x},\bar{u}\right)  \in\mbox{\rm Graph}\left(  S\right)  $.
Then $S$ has the Aubin property at $\bar{x}$ with respect to
$\bar{u}$ if and only if $D^{\ast
}S(\bar{x}\mid\bar{u})(\mathbf{0})=\{\mathbf{0}\}$ or equivalently
$|D^{\ast}S(\bar{x}\mid\bar {u})|^{+}<\infty$. In this case,
$\mbox{\rm lip}\,S\left(  \bar{x}\mid \bar{u}\right)
=|D^{\ast}S(\bar{x}\mid\bar{u})|^{+}$.
\end{proposition}

Using the above criterion we show that an everywhere continuous
strictly increasing single-valued map from the reals to the reals
could be nowhere Lipschitz continuous.

\begin{example}
Let $A\subset\mathbb{R}$ be a measurable set with the property that
for every $a,b\in\mathbb{R}$, $a<b$, the Lebesgue measure of
$A\cap\left(  a,b\right)  $ satisfies $0<m(A\cap\lbrack
a,b])<\left\vert b-a\right\vert $. Consider the function $f:\left[
0,1\right]  \rightarrow\mathbb{R}$ defined by $f\left( x\right)
=m\left(  A\cap\left(  0,x\right)  \right)  .$ Note that the
derivative $f^{\prime}\left(  x\right)  $ exists almost everywhere
and is equal to $\chi_{A}\left(  x\right)  $, the characteristic
function of $A$ (equal to $1$ if $x\in A$ and $0$ if not). This
means that every point $\bar {x}\in\left[  0,1\right]  $ is
arbitrarily close to a point $x$ where $f^{\prime}\left(  x\right) $
is well-defined and equals zero. Thus $(0,1)\in
N_{\mathrm{Graph}(f)}(\bar{x},f(\bar{x}))$. The function $f$ is
strictly increasing and continuous, so it has a continuous inverse
$g:\left[ 0,f\left( 1\right)  \right]  \rightarrow\left[  0,1\right]
$. Applying the Mordukhovich criterion
(Proposition~\ref{pro:Boris-criterion}) we obtain that $g$ does not
have the Aubin property at $f\left(  \bar{x}\right)  $. It follows
that $g$ is not strictly continuous at $f\left(  \bar{x}\right)  $
and in fact neither is so at any $y\in\left[  0,f\left(  1\right)
\right]  $.$\hfill\Box$
\end{example}

\section{Preliminary results in Variational Analysis}
\label{sec-3}

In this section we establish a Sard type result for the image of the
set of local minima (respectively, local Pareto minima) in case of
single--valued scalar (respectively, vector--valued) functions. We
also obtain several auxiliary results that will be used in
Section~\ref{sec-main}.

\subsection{Sard result for local (Pareto) minima}

In this subsection we use simple properties on the continuity of
set-valued maps to obtain a Sard type result for local minima for
both scalar and vector-valued functions. Let us recall that a
(single-valued) function $f:X\rightarrow\mathbb{R}$ is called
\emph{lower semicontinuous} \emph{at $\bar{x}$} if
\[
\underset{x\rightarrow\bar{x}}{\lim\inf}\,f(x)\geq f(\bar{x})\,.
\]
The function $f$ is called \emph{lower semicontinuous}, if it is
lower semicontinuous at every $x\in X.$ It is well-known that a
function $f$ is lower semicontinuous if and only if its sublevel
sets $$[f\leq r]:=\{x\in X:f(x)\leq r\}$$ are closed for all
$r\in\mathbb{R}$.

\begin{proposition}
[Sublevel map]\label{pro:level-set-cts} Let $D$ be a closed subset
of a complete metric space $X$ and $f:D\rightarrow\mathbb{R}$ be a
lower semicontinuous function. Then the (sublevel) set-valued map
\[
\left\{
\begin{array}
[c]{l}
L_{f}:\mathbb{R}\rightrightarrows D\smallskip\\
L_{f}\left(  r\right)  =[f\leq r]\cup\partial D
\end{array}
\right.
\]
is outer semicontinuous. Moreover, $L_{f}$ is continuous at
$\bar{r}\in f\left( D\right)$ if and only if there is no
$x\in\mbox{\rm int}\left( D\right)$ such that $f\left(  x\right)
=\bar{r}$ and $x$ is a local minimizer of $f$.
\end{proposition}

\noindent\textit{Proof.} The map $L_{f}^{\prime}:\mathbb{R}\rightrightarrows
D$ defined by $L_{f}^{\prime}(r)=f^{-1}\left(  (-\infty,r]\right)  $ is outer
semicontinuous since $f$ is lower semicontinuous (see \cite[Example 5.5]{RW98}
for example), so $L_{f}$ is easily seen to be outer semicontinuous.

We now prove that $L_{f}$ is inner semicontinuous at $\bar{r}$ under
the additional conditions mentioned in the statement. For any
$r_{i}\rightarrow\bar{r}$, we want to show that if $\bar{x}\in
L_{f}\left(  \bar{r}\right)  $, then there exists
$x_{i}\rightarrow\bar{x}$ such that $x_{i}\in L_{f}\left(
r_{i}\right)  $. We can assume that $f\left(  \bar{x}\right)
=\bar{r}$ and $r_{i}<\bar{r}$ for all $i$, otherwise we can take
$x_{i}=\bar{x}$ for $i$ large enough. Since $\bar{x}$ is not a local
minimum, for any $\epsilon>0$, there exists $\delta>0$ such that if
$\left\vert \bar{r}-r_{i}\right\vert <\delta$, there exists an
$x_{i}$ such that $f\left(  x_{i}\right)  \leq r_{i}$ and
$\left\vert x_{i}-\bar{x}\right\vert <\epsilon$.

For the converse, assume now that $L_{f}$ is inner semicontinuous at
$\bar{r}$. Then taking $r_{i}\nearrow\bar{r}$ we obtain that for
every $x\in\mbox{int}\left( D\right)  \cap f^{-1}\left(
\bar{r}\right)  $, there exists $x_{i}\in f^{-1}\left(  r_{i}\right)
$ with $x_{i}\rightarrow x$. Since $f\left( x_{i}\right)
=r_{i}<\bar{r}=f\left(  x\right)  $, $x$ cannot be a local minimum.
$\hfill\Box$

\bigskip

According to the above result, if $f$ has no local minima, then the set-valued
map $L_{f}$ is continuous everywhere. The above result has the following
interesting consequence.

\begin{corollary}
[Local minimum values]\label{corollary:lsc-generic-min-vals}Let $M_{f}$ denote
the set of local minima of a lower semicontinuous function $f:D\rightarrow
\mathbb{R}$ (where $D$ is a closed subset of a complete space $X$). Then the
set $f(M_{f}\cap\mbox{int}\left(  D\right)  )$ is meager in $\mathbb{R}$.
\end{corollary}

\noindent\textit{Proof.} Since the set-valued map $L_{f}$ (defined
in Proposition~\ref{pro:level-set-cts}) is outer semicontinuous
(with closed-values), it is generically continuous by
Theorem~\ref{theorem:[RW98,5.55]}. The second part of Proposition
\ref{pro:level-set-cts} yields the result on $f$. $\hfill\Box$

\bigskip

It is interesting to compare the above result with the classical
Sard theorem. We recall that the Sard theorem asserts that the image
of critical points (derivative not surjective) of a $C^{k}$ function
$f:\mathbb{R}^{n} \rightarrow\mathbb{R}^{m}$ is of measure zero
provided $k>n-m.$ (See \cite{Sard42}; the case $m=1$ is known as the
Sard-Brown theorem \cite{Brown35}.)
Corollary~\ref{corollary:lsc-generic-min-vals} asserts the
topological sparsity of the (smaller) set of minimum values for
scalar functions ($m=1$), without assuming anything but lower
semicontinuity (and completeness of the domain).

\bigskip

We shall now extend Corollary~\ref{corollary:lsc-generic-min-vals}
in the vectorial case. We recall that a set $K\subset\mathbb{R}^{m}$
is a \emph{cone}, if $\lambda K\subset K$ for all $\lambda\geq0.$ A
cone $K$ is called \emph{pointed} if $K\cap(-K)=\{\mathbf{0}_m\}$
(or equivalently, if $K$ contains no full lines). It is well-known
that there is a one-to-one correspondence between pointed convex
cones of $\mathbb{R}^{m}$ and partial orderings in $\mathbb{R}^{m}$.
In particular, given such a cone $K$ of $\mathbb{R}^{m}$ we set
$y_{1}\leq _{K}y_{2}$ if and only if $y_{2}-y_{1}\in K$ (see for
example, \cite[Section~3E]{RW98}). Further, given a set-valued map
$S:\mathbb{R}^{n}\rightrightarrows\mathbb{R}^{m}$ we say that

\begin{itemize}
\item $\bar{x}$ is a \emph{(local) Pareto minimum} of $S$ with \emph{(local)
Pareto minimum value} $\bar{y}$ if there is a neighborhood $U$ of
$\bar{x}$ such that if $x\in U$ and $y\in S\left(  x\right)  $, then
$y\not \leq_{K}\bar{y}$, i.e., $S\left(  U\right)  \cap\left(
\bar{y}-K\right)  =\bar {y}$.
\end{itemize}

For $S:\mathbb{R}^{n}\rightrightarrows\mathbb{R}^{m}$, define the
map $S_{K}:\mathbb{R}^{n}\rightrightarrows\mathbb{R}^{m}$ by
$S_{K}\left(x\right)=S\left(x\right)+K$. The graph of $S_{K}$ is
also known as the \emph{epigraph } \cite{GRTZ03,Jahn04} of $S$. One
easily checks that $y\in S_{K}\left(x\right)$ implies $y+K\subset
S_{K}\left(x\right)$. Here is our result on local Pareto minimum
values.

\begin{proposition}
[Pareto minimum values] Let
$S:\mathbb{R}^{n}\rightrightarrows\mathbb{R}^{m}$ be an outer
semicontinuous map such that $y\in S\left( x\right)$ implies
$y+K\subset S\left(  x\right)$ (that is, $S=S_K$). Then the set of
local Pareto minimum values is meager.
\end{proposition}

\noindent\textit{Proof.} Since $S$ is outer semicontinuous, then
$S^{-1}$ is outer semicontinuous as well by
\cite[Theorem~5.7(a)]{RW98}, so $S^{-1}$ is generically continuous
by Theorem~\ref{theorem:[RW98,5.55]}. Suppose that $\bar{y}$ is a
local Pareto minimum of a local Pareto minimizer $\bar{x}$.

By the definition of local Pareto minimum, there is a neighborhood
$U$ of $\bar{x}$ such that if $y\leq_{K}\bar{y}$ and $y\neq\bar{y}$,
then $S^{-1}\left( y\right) \cap U=\emptyset$. (We can assume that
$y$ is arbitrarily close to $\bar{y}$ since $S^{-1}\left(  y\right)
\subset S^{-1}\left(  \lambda y+\left(  1-\lambda\right)
\bar{y}\right)  $ for all $0\leq\lambda\leq1$.) Therefore,
$\bar{x}\notin\liminf_{y\rightarrow\bar{y}}S^{-1}\left(  y\right) $.
In other words, $S^{-1}$ is not continuous at $\bar{y}$. Therefore,
the set of local Pareto minimum values is meager. $\hfill\Box$

\bigskip

We show how the above result compares to critical point results. Let
us recall from \cite{Ioffe08} the definition of critical points of a
set-valued map. Given a metric space $X$ (equipped with a distance
$\rho$) we denote by $B_{\rho}(x,\lambda)$ the set of all $x'\in X$
such that $\rho(x,x')\le \lambda$.

\begin{definition}\label{Def-critical}
Let $(X,\rho_1)$ and $(Y,\rho_2)$ be metric spaces, and let
$S:X\rightrightarrows Y$. For $\left(x,y\right)\in\gph\,(S)$, we set
\[
\Sur\, S\left(x\mid y\right)\left(\lambda\right)=\sup\left\{
r\geq0\mid B_{\rho_{2}}(y,r)\subset
S\left(B_{\rho_{1}}(x,\lambda\right)\right\}
\]
and then for $\left(\bar{x},\bar{y}\right)\in\gph\,(S)$ define the
\emph{rate of surjection }of $S$ at $\left(\bar{x},\bar{y}\right)$
by
\[
\sur\,
S\left(\bar{x}\mid\bar{y}\right)=\liminf_{\left(x,y,\lambda\right)\rightarrow\left(\bar{x},\bar{y},+0\right)}\frac{1}{\lambda}\,\Sur\,
S\left(x\mid y\right)\left(\lambda\right).
\]
We say that $S$ is \emph{critical} at $\left(\bar x,\bar
y\right)\in\gph\,(S)$ if $\sur S\,\left(\bar x\mid \bar y\right)=0$,
and regular otherwise. Also, $\bar y$ is a \emph{(proper) critical
value} of $S$ if there exists $\bar x$ such that $\bar y\in
S\left(\bar x\right)$ and $S$ is critical at $\left(\bar x,\bar
y\right)$.
\end{definition}

This definition of critical values characterizes the values at which
metric regularity is absent. In the particular case where
$S:\mathbb{R}^n\rightarrow\mathbb{R}^m$ is a $\mathcal{C}^1$
function, critical points correspond exactly to where the Jacobian
has rank less than $m$. We refer to \cite{Ioffe08} for more details.

\smallskip

One easily sees that if $y$ is a Pareto minimum value of $S$, then
there exists $x\in X$ such that $\left(x,y\right)\in\gph\,(S)$, and
$\Sur\,S\left(x\mid y\right)\left(\lambda\right)=0$ for all small
$\lambda>0$. This readily implies that $y$ is a critical value.

\subsection{Extending the Mordukhovich criterion and a critical value result}

The two results of this subsection are important ingredients of the
forthcoming proof of our main theorem. The first result we need is
an adaptation of the Mordukhovich criterion
(Proposition~\ref{pro:Boris-criterion}) to the case where the domain
of a set-valued function $S$ is (included in) a smooth submanifold
$\mathcal{X}$ of $\mathbb{R}^{n}$. (Note that this new statement
recovers the Mordukhovich criterion if
$\mathcal{X}=\mathbb{R}^{n}$.)

\begin{proposition}
\label{pro:Mordukhovich-extended}(Extended Mordukhovich criterion)
Let $\mathcal{X}\subset\mathbb{R}^{n}$ be a $\mathcal{C}^{1}$ smooth
submanifold of dimension $d$ and
$S:\mathcal{X}\rightrightarrows\mathbb{R}^{m}$ be a set-valued map
whose graph is locally closed at $\left( \bar{x},\bar{y}\right)  \in
\mbox{\rm Graph}\left(  S\right)  $. Consider the mapping
\[
\left\{
\begin{array}
[c]{l}
H:\mathbb{R}^{m}\rightrightarrows\mathbb{R}^{n}\smallskip\\
H\left(  y^{\ast}\right)  =D^{\ast}S\left( \bar{x}\mid\bar{y}\right)
\left( y^{\ast}\right) \, \cap
\,T_{\mathcal{X}}\left(\bar{x}\right).
\end{array}
\right.
\]
If $H\left(  \mathbf{0}_m\right)  =\left\{  \mathbf{0}_n\right\}  $,
or equivalently
\[
N_{{\scriptsize \mbox{\rm Graph}\left(  S\right)  }}\left(
\bar{x},\bar {y}\right)  \cap\left(  T_{\mathcal{X}}\left(
\bar{x}\right)  \times\left\{ \mathbf{0}_{m}\right\}  \right)
=\left\{ \mathbf{0}_{n+m}\right\}  ,
\]
then $S$ has the Aubin property at $\bar{x}$ for $\bar{y}$ relative to
$\mathcal{X}$. Furthermore,
\[
\mbox{\rm lip}_{\mathcal{X}}S\left(  \bar{x}\mid\bar{y}\right)
=\left\vert H\right\vert ^{+}=\sup\left\{  \frac{\left\vert
u\right\vert }{\left\vert v\right\vert }\mid\left(  u,v\right)  \in
N_{{\scriptsize \mbox{\rm Graph}\left(  S\right)  }}\left(
\bar{x},\bar {y}\right)  \cap\left(  T_{\mathcal{X}}\left(
\bar{x}\right)  \times \mathbb{R}^{m}\right)  \right\}  .
\]
\end{proposition}

\noindent\textit{Proof.} Fix $\left(  \bar{x},\bar{y}\right)  \in
\mbox{\rm Graph}\left(  S\right)  $ and denote by $N_{\mathcal{X}}\left(
\bar{x}\right)  $ the normal space of $\mathcal{X}$ at $\bar{x}$ (seing as
subspace of $\mathbb{R}^{n},$ that is, $T_{\mathcal{X}}\left(  \bar{x}\right)
\oplus N_{\mathcal{X}}\left(  \bar{x}\right)  =\mathbb{R}^{n}$). Given a
closed neighborhood $U$ of $\left(  \bar{x},\bar{y}\right)  $, we define the
function
\[
\left\{
\begin{array}
[c]{l}
{\tilde{S}}:\mathbb{R}^{n}\rightrightarrows\mathbb{R}^{m}\smallskip\\
\mbox{\rm Graph}\left(  {\tilde{S}}\right)  =\left(  \mbox{\rm
Graph}\left( S\right)  \cap U\right)  +\left(  N_{\mathcal{X}}\left(
\bar{x}\right) \times\left\{  \mathbf{0}_{m}\right\}  \right)  .
\end{array}
\right.
\]
Shrinking the neighborhood $U$ around $\left( \bar{x},\bar{y}\right)
$ if necessary, we may assume that every $\left(  x,y\right)  \in U$
can be represented uniquely as a sum of elements in $\left(
\mathcal{X} \times\mathbb{R}^{m}\right)  \cap U$ and
$N_{\mathcal{X}}\left(  \bar {x}\right)  \times\left\{
\mathbf{0}_{m}\right\}  $. Since $\mbox{\rm Graph}\left(  S\right) $
is locally closed, we can choose $U$ small enough so that $\mbox{\rm
Graph}\left(  S\right)  \cap U$ is closed. Further, since $\mbox{\rm
Graph}\left(  \tilde{S}\right)  $ is homeomorphic to $\left(
\mbox{\rm Graph}\left(  S\right)  \cap U\right)  \times
\mathbb{R}^{n-d}$, it is also closed.

\medskip

\textbf{Step 1: (Relating $\tilde{S}$ to $H$)} By applying a result
on the normal cones under set addition \cite[Exercise~6.44]{RW98},
we have $N_{{\scriptsize \mbox{\rm Graph}\left(  \tilde{S}\right)
}}\left(  \bar {x},\bar{y}\right)  \subset N_{{\scriptsize \mbox{\rm
Graph}\left(  S\right) }}\left(  \bar{x},\bar{y}\right)  \cap\left(
T_{\mathcal{X}}\left(  \bar {x}\right) \times\mathbb{R}^{m}\right)$.
To prove the reverse inclusion, note that for every $\left(
x,y\right)  \in\mbox{\rm Graph}\left( \tilde {S}\right)$ near
$\left(\bar{x},\bar{y}\right)$ with $\left(x,y\right)  =\left(
x_{1},y\right)+\left(x_{2},\mathbf{0}_{m}\right)$, where $\left(
x_{1},y\right) \in\mbox{\rm Graph}\left(  S\right) $ and $x_{2}\in
N_{\mathcal{X}}\left(  \bar{x}\right)$, one easily sees that
$\hat{N}_{{\scriptsize \mbox{\rm Graph}\left( \tilde{S}\right)
}}\left( x,y\right)  \supset\hat{N}_{{\scriptsize \mbox{\rm
Graph}\left(  S\right)  } }\left(  x_{1},y\right) \cap\left(
T_{\mathcal{X}}\left(  \bar{x}\right) \times\mathbb{R}^{m}\right)$.
The extension of this inclusion to limiting normal cones is
immediate. Therefore we obtain
\[
N_{{\scriptsize \mbox{\rm Graph}\left(  \tilde{S}\right)  }}\left(
\bar {x},\bar{y}\right)  =N_{{\scriptsize \mbox{\rm Graph}\left(
S\right)  } }\left(  \bar{x},\bar{y}\right)  \cap\left(
T_{\mathcal{X}}\left(  \bar {x}\right) \times\mathbb{R}^{m}\right),
\]
and so $D^{\ast}\tilde{S}\left(  \bar{x}\mid\bar{y}\right)$ equals
the set-valued map $H$ described in the statement. Thus
\begin{eqnarray*}
D^{*}\tilde{S}\left(\bar{x}\mid\bar{y}\right)\left(\mathbf{0}_m\right)
& = & \left\{ x^{*}\mid\left(x^{*},\mathbf{0}_m\right)\in N_{\scriptsize\gph\left(\tilde{S}\right)}\left(\bar{x},\bar{y}\right)\right\} \\
& = & \left\{ x^{*}\mid\left(x^{*},\mathbf{0}_m\right)\in N_{\scriptsize\gph\left(S\right)}\left(\bar{x},\bar{y}\right)\cap\left(T_{\mathcal{X}}
\left(\bar{x}\right)\times\mathbb{R}^{m}\right)\right\} \\
& = & \left\{ \mathbf{0}_n\right\} ,
\end{eqnarray*}
and by the Mordukhovich criterion, the map $\tilde{S}$ has the Aubin
property at $\bar{x}$ for $\bar{y}$. \smallskip

Taking neighborhoods $V$ of $\bar{x}$ and $W$ of $\bar{y}$ so that $S\left(
x\right)  \cap W=\tilde{S}\left(  x\right)  \cap W$ for all $x\in
V\cap\mathcal{X}$, we deduce that $S$ has the Aubin property at $\bar{x}$ for
$\bar{y}$ relative to $\mathcal{X}$ as asserted.

\medskip

\textbf{Step 2: ($\mbox{\rm lip}_{\mathcal{X}}S\left(
\bar{x}\mid\bar {y}\right)  =\left\vert H\right\vert ^{+}$)} The
Mordukhovich criterion on $\tilde{S}$ yields $$\left\vert
H\right\vert ^{+}=\mbox{\rm lip}\,\tilde {S}\left(
\bar{x}\mid\bar{y}\right)  \geq\mbox{\rm lip}_{\mathcal{X}
}S\left(  \bar{x}\mid\bar{y}\right)\,.$$ Our task is thus to prove
that the above inequality is actually an equality. Since $\mbox{\rm
lip}\,\tilde{S}\left( \bar{x}\mid\bar{y}\right)  =\left\vert
H\right\vert ^{+}$, for any $\kappa<\left\vert H\right\vert ^{+}$
and neighborhoods $V$ of $\bar{x}$ and $W$ of $\bar{y}$, there exist
$x_{1},x_{2}\in V$ such that
\[
\tilde{S}\left(  x_{2}\right)  \cap W\not \subset \tilde{S}\left(
x_{1}\right)  +\kappa\left\vert x_{1}-x_{2}\right\vert \mathbb{B}.
\]
Note that $\tilde{S}\left(  x_{1}\right)  =\tilde{S}\left(  P\left(
x_{1}\right)  \right)  $, $\tilde{S}\left(  x_{2}\right)
=\tilde{S}\left( P\left(  x_{2}\right)  \right)  $ and $\left\vert
P\left(  x_{1}\right) -P\left(  x_{2}\right)  \right\vert
\leq\left\vert x_{1}-x_{2}\right\vert $, where $P$ stands for the
projection of $\mathbb{R}^{n}$ onto $\bar {x}+T_{\mathcal{X}}\left(
\bar{x}\right)$. We may choose $V$ to be a ball containing
$\bar{x}$, and define the projection parametrization $L:\left(
\bar{x}+T_{\mathcal{X}}\left( \bar{x}\right) \right) \cap
V\rightarrow\mathcal{X}$ of the manifold $\mathcal{X}$ by the
relation $x-L\left(  x\right)  \in N_{\mathcal{X}}\left(
\bar{x}\right)$. Shrinking $V$ if needed, the map $L$ becomes
single-valued and smooth (in fact, it is a local chart of
$\mathcal{X}$ at $\bar{x}$ provided we identify $\bar
{x}+T_{\mathcal{X}}\left(  \bar{x}\right)  $ with $\mathbb{R}^{d}$).
Furthermore, $L$ has Lipschitz constant $1$ at $\bar{x}$. Therefore,
for any $\epsilon>0$, we can reduce $V$ as needed so that $L$ is
Lipschitz continuous in its domain with Lipschitz constant at most
$\left( 1+\epsilon\right)  $ using standard arguments ({\it e.g.}
\cite[Thms~9.7,~9.2]{RW98}). This means that
\[
S\left(  L\left(  x_{2}\right)  \right)  \cap W=\tilde{S}\left(  x_{2}\right)
\cap W\not \subset \tilde{S}\left(  x_{1}\right)  +\kappa\left\vert
x_{1}-x_{2}\right\vert \mathbb{B}=S\left(  L\left(  x_{1}\right)  \right)
+\kappa\left\vert x_{1}-x_{2}\right\vert \mathbb{B}.
\]
By the Lipschitz continuity of $L$, we have $\left\vert L\left(
x_{1}\right) -L\left(  x_{2}\right)  \right\vert \leq\left(
1+\epsilon\right)  \left\vert x_{1}-x_{2}\right\vert $, which gives
\[
S\left(  L\left(  x_{2}\right)  \right)  \cap W\not \subset S\left(  L\left(
x_{2}\right)  \right)  +\frac{\kappa}{\left(  1+\epsilon\right)  }\left\vert
L\left(  x_{1}\right)  -L\left(  x_{2}\right)  \right\vert \mathbb{B},
\]
yielding
\[
\frac{\kappa}{1+\epsilon}\leq\mbox{\rm lip}_{\mathcal{X}}S\left(  \bar
{x}\mid\bar{y}\right)  .
\]
Since $\kappa$ and $\epsilon$ are arbitrary, we conclude that
$\left\vert H\right\vert ^{+}=\mbox{\rm lip}_{\mathcal{X}}S\left(
\bar{x}\mid\bar {y}\right)$ as asserted.

\smallskip

The proof is complete.$\hfill\Box$

\medskip
The second result is an adaptation of part of
\cite[Theorem~6]{Ioffe08}. Recall that for a smooth function
$f:\mathbb{R}^{n}\rightarrow\mathbb{R}^{m}$,
$\bar{x}\in\mathbb{R}^{n}$ is a \emph{critical point} if the
derivative $\nabla f(\bar{x})$ is not surjective, while $\bar{y}
\in\mathbb{R}^{m}$ is a \emph{critical value} if there is a critical
point $\bar{x}$ for which $f\left( \bar{x}\right)=\bar{y}$. (Note
this is a particular case of the general definition given in
Definition~\ref{Def-critical}.)

\begin{lemma}
\label{lemma:Ioffes-lemma}Let $\mathcal{X}$ be a $\mathcal{C}^{k}$
smooth manifold in $\mathbb{R}^{n}$ of dimension $d$, and
$\mathcal{M}$ be a $\mathcal{C}^{k}$ manifold in $\mathbb{R}^{n+m}$
such that $\mathcal{M}\subset\mathcal{X} \times\mathbb{R}^{m}$, with
$k>\dim\mathcal{M}-\dim\mathcal{X}$. Then the set of points
$x\in\mathcal{X}$ such that there exists some $y$ satisfying $\left(
x,y\right)  \in\mathcal{M}$ and $N_{\mathcal{M}}\left(  x,y\right)
\cap\left(  T_{\mathcal{X}}\left(  x\right)  \times\left\{
\mathbf{0}_{m} \right\}  \right)  \supsetneq\left\{
\mathbf{0}_{n+m}\right\} $ is of Lebesgue measure zero
in~$\mathcal{X}$.
\end{lemma}

\noindent\textit{Proof.} Let $\mathrm{Proj}_{\mathcal{M}}$ denote
the restriction to the manifold $\mathcal{M}$ of the projection of
$\mathcal{X} \times\mathbb{R}^{m}$ onto $\mathcal{X}$. As
$k>\dim\mathcal{M}-\dim \mathcal{X}$, the set of critical values of
$\mathrm{Proj}_{\mathcal{M}}$ is a set of measure zero by the
classical Sard theorem \cite{Sard42}. Let $\left( x,y\right)
\in\mathcal{M}$ and assume $\left(  x^{\ast},\mathbf{0}_m\right) \in
N_{\mathcal{M}}\left(  x,y\right)  \cap\left(  T_{\mathcal{X}}\left(
x\right)  \times\left\{  \mathbf{0}_{m}\right\}  \right)$ with
$x^{\ast} \neq\mathbf{0}_n$. This gives
\[
T_{\mathcal{M}}\left(  x,y\right)  =\left(  N_{\mathcal{M}}\left(  x,y\right)
\right)  ^{\perp}\subset\left\{  x^{\ast}\right\}  ^{\perp}\times
\mathbb{R}^{m},
\]
where $\left\{  x^{\ast}\right\}  ^{\perp}=\left\{  x^{\prime}\in
\mathbb{R}^{n}\mid\left\langle x^{\ast},x^{\prime}\right\rangle
=0\right\}  $. Since $T_{\mathcal{M}}\left(  x,y\right)  \subset
T_{\mathcal{X}}\left( x\right)  \times\mathbb{R}^{m}$ we obtain
\[
T_{\mathcal{M}}\left(  x,y\right)  \subset\left(  \left\{  x^{\ast}\right\}
^{\perp}\cap T_{\mathcal{X}}\left(  x\right)  \right)  \times\mathbb{R}^{m}.
\]
Let $Z$ stand for the subspace on the right hand side. Then the
projection of $Z$ onto $T_{\mathcal{X}}\left(  x\right)  $ is a
proper subspace of $T_{\mathcal{X}}\left(  x\right)  $. All the
more, this applies to $T_{\mathcal{M}}\left(  x,y\right)  $. By
\cite[Corollary 3]{Ioffe08}, this implies that $\left(  x,y\right)
$ is a singular point of $\mathrm{Proj} _{\mathcal{M}}$, so $x$ is a
critical value of $\mathrm{Proj}_{\mathcal{M}}$. The conclusion of
the lemma follows. $\hfill\Box$

\medskip

\subsection{Linking sets}

We introduce the notion of \emph{linking} that is commonly used in
critical point theory. Let us fix some terminology: if
$B\subset\mathbb{R}^{n}$ is homeomorphic to a subset of
$\mathbb{R}^{d}$ with nonempty interior, we say that the set
$\partial B$ is the \emph{relative boundary} of $B$ if it is a
homeomorphic image of the boundary of the set in $\mathbb{R}^{d}$.

\begin{definition}
\label{Def:link}\cite[Section~II.8]{Struwe00} Let $A$ be a subset of
$\mathbb{R}^{n+m}$ and let $B$ be a submanifold of
$\mathbb{R}^{n+m}$ with relative boundary $\partial B$. Then we say
that $A$ and $\Gamma=\partial B$ \emph{link }if \smallskip

(i) $A\cap\Gamma=\emptyset$ \smallskip

(ii) for any continuous map $h\in\mathcal{C}^{0}\left(  \mathbb{R}
^{n+m},\mathbb{R}^{n+m}\right)  $ such that $h\mid_{\Gamma}=id$ we
have $h\left(  B\right)  \cap A\neq\emptyset$.
\end{definition}

In particular, the following result holds. This result will be used in
Section~\ref{sec-main}.

\begin{theorem}
[Linking sets]\label{theorem:linking-sets} Let $\mathcal{K}_{1}$ and
$\mathcal{K}_{2}$ be linear subspaces such that $\mathcal{K}_{1}
\oplus\mathcal{K}_{2}=\mathbb{R}^{n+m}$, and take any $\bar{v}\in
\mathcal{K}_{1}\backslash\left\{ \mathbf{0}\right\}$. Then for
$0<r<R$, the sets
$$
A:=\mathbb{S}\left( \mathbf{0},r\right) \cap\mathcal{K}_{1} \qquad
\text{and}\qquad \Gamma:=\left( \mathbb{B}\left( \mathbf{0},R\right)
\cap \mathcal{K}_{2}\right) \cup\left( \mathbb{S}\left(
\mathbf{0},R\right) \cap\left( \mathcal{K}_{2}+\mathbb{R}_{+}\left\{
\bar{v}\right\}  \right) \right) $$ link.
\end{theorem}

\noindent\textit{Proof.} Use methods in \cite[Section II.8]{Struwe00}, or
infer from Example 3 there. $\hfill\Box$

\bigskip

We finish this section with two useful results. The first one is well-known
(with elementary proof) and is mentioned for completeness.

\begin{proposition}
\label{pro:orthogonal-intersection-sum} If $\mathcal{K}_{1}$ and
$\mathcal{K}_{2}$ are subspaces of $\mathbb{R}^{n+m}$, then
$\mathcal{K} _{1}^{\perp}\cap\mathcal{K}_{2}^{\perp}=\left\{
\mathbf{0}\right\}  $ if and only if
$\mathcal{K}_{1}+\mathcal{K}_{2}=\mathbb{R}^{n+m}$.
\end{proposition}

The following lemma will be needed in the proof of forthcoming
Lemma~\ref{lemma:normals-behave} (Section~\ref{sec-main}).

\begin{lemma}
\label{lemma:extension-to-interior} If the sets $\mathbb{B}\left(
\mathbf{0},1\right)  $ and $D$ are homeomorphic, then any
homeomorphism $f$ between $\mathbb{S}\left(\mathbf{0},1\right)$ and
$\partial D$ can be extended to a homeomorphism $F:\mathbb{B}\left(
\mathbf{0},1\right) \rightarrow D$ so that $F\mid_{\mathbb{S}\left(
\mathbf{0},1\right) }=f$.
\end{lemma}

\noindent\textit{Proof.} Let $H:\mathbb{B}\left( \mathbf{0},1\right)
\rightarrow D$ be a homeomorphism between $\mathbb{B}\left(
\mathbf{0} ,1\right)  $ and $D$ and denote $h:\mathbb{S}\left(
\mathbf{0},1\right) \rightarrow\partial D$ by
$h=H\mid_{\mathbb{S}\left(  \mathbf{0},1\right)  }$. We define the
(continuous) function $F:\mathbb{B}\left(  \mathbf{0},1\right)
\rightarrow D$ by
\[
F\left(  x\right)  =
\begin{cases}
H\left(  \left\vert x\right\vert h^{-1}(f\left(  x/|x|\right)  \right)
\smallskip & \mbox{if }x\neq\mathbf{0}\\
\qquad H\left(  \mathbf{0}\right)  & \mbox{if }x=\mathbf{0}.
\end{cases}
\]
It is straightforward to check that $F\mid_{\mathbb{S}\left(
\mathbf{0} ,1\right)  }=f$. Let us show that $F$ is injective:
indeed, if $F\left( x_{1}\right)  =F\left(  x_{2}\right)  $, then
$\left\vert x_{1}\right\vert h^{-1}(f\left(  x_{1}/|x_{1}|\right)
)=\left\vert x_{2}\right\vert h^{-1}(f\left(  x_{2}/|x_{2}|\right)
).$ If both sides are zero, then $x_{1}=x_{2}=\mathbf{0}$. Otherwise
$\left\vert x_{1}\right\vert =\left\vert x_{2}\right\vert $ and
$x_{1}/|x_{1}|=x_{2}/|x_{2}|$, which implies that
$x_{1}=x_{2}$.\smallskip\newline To see that $F$ is a bijection, fix
any $y\in D$, and let $x^{\prime}\in\mathbb{B}\left(
\mathbf{0},1\right)  $ be such that $y=H\left(  x^{\prime}\right) $.
If $x^{\prime}=\mathbf{0}$, then $y=F\left(  \mathbf{0}\right)  $.
Otherwise,
\[
y=H\left(  \left\vert x^{\prime}\right\vert \left(  \frac{x^{\prime}
}{|x^{\prime}|}\right)  \right)  =H(\left\vert x^{\prime}\right\vert
\,h^{-1}\circ f\left(  f^{-1}\circ
h(\frac{x^{\prime}}{|x^{\prime}|})\right) =F\left(  \left\vert
x^{\prime}\right\vert f^{-1}\circ h(\frac{x^{\prime}
}{|x^{\prime}|})\right)  .
\]
This shows that $F$ is also surjective, thus a continuous bijection. Since
$\mathbb{B}\left(  \mathbf{0},1\right)  $ is compact, it follows that $F$ is a
homeomorphism.$\hfill\Box$

\section{Generic continuity of tame set-valued maps}

\label{sec-main}

From now on we limit our attention to the class of semialgebraic (or
more generally, o-minimal) set-valued maps. In this setting our main
result eventually asserts that every such set-valued map is
generically strictly continuous (see Section~\ref{subsec-3}). To
prove this, we shall need several technical lemmas, given in
Section~\ref{subsec-2}. In Section~\ref{subsec-1} we give
preliminary definitions and results of our setting.

\subsection{Semialgebraic and definable mappings}

\label{subsec-1}

In this section we recall basic notions from semialgebraic and o-minimal
geometry. Let us define properly the notion of a semialgebraic set
(\cite{BR90}, \cite{Coste99}). (We denote by $\mathbb{R}[x_{1},\ldots,x_{n}]$
the ring of real polynomials of $n$ variables.)

\begin{definition}
[Semialgebraic set]\label{Def:semialgebraic}A subset $A$ of
$\mathbb{R}^{n}$ is called \emph{semialgebraic} if it has the form
\[
A=\bigcup\limits_{i=1}^{k}\{x\in\mathbb{R}^{n}:p_{i}(x)=0,q_{i1}
(x)>0,\ldots,q_{i\ell}(x)>0\},
\]
where $p_{i},q_{ij}\in\mathbb{R}[x_{1},\ldots,x_{n}]$ for all $i\in
\{1,\ldots,k\}$ and $j\in\{1,\ldots,\ell\}$.
\end{definition}

In other words, a set is semialgebraic if it is a finite union of
sets that are defined by means of a finite number of polynomial
equalities and inequalities. A set-valued map
$S:\mathbb{R}^{n}\rightrightarrows \mathbb{R}^{m}$ is called
\emph{semialgebraic}, if its graph $\mbox{\rm Graph}\left( S\right)
$ is a semialgebraic subset of $\mathbb{R}^{n}\times\mathbb{R}^{m}
$.\smallskip

An important property of semialgebraic sets is that of Whitney
stratification (\cite[§4.2]{vdDries98},
\cite[Theorem~6.6]{Coste99}).

\begin{theorem}
\label{theorem:stratification} ($\mathcal{C}^{k}$ stratification)
For any $k\in\mathbb{N}$ and any semialgebraic subsets
$X_{1},\dots,X_{l}$ of $\mathbb{R}^{n}$, we can write
$\mathbb{R}^{n}$ as a disjoint union of finitely many semialgebraic
$\mathcal{C}^{k}$ manifolds $\{\mathcal{M} _{i}\}_{i}$ (that is,
$\mathbb{R}^{n}=\dot{\cup}_{i=1}^{I}\mathcal{M}_{i}$) so that each
$X_{j}$ is a finite union of some of the $\mathcal{M}_{i}$'s.
Moreover, the induced stratification $\{\mathcal{M}_{i}^{j}\}_{i}$
of $X_{j}$ has the Whitney property, that is, for any sequence
$\{x_{\nu}\}_{\nu} \subset\mathcal{M}_{i}^{j}$ converging to
$x\in\mathcal{M}_{i_{0}}^{j}$ we have
$$\underset{v\rightarrow\infty}{\lim\sup\,}N_{\mathcal{M}_{i}^{j}}(x_{\nu
})\subset N_{\mathcal{M}_{i_{0}}^{j}}(x_{\nu}).$$
\end{theorem}

\noindent In particular, every semi-algebraic set can be written as
a finite disjoint union of manifolds (\textquotedblleft
strata\textquotedblright) that fit together in a regular way
(\textquotedblleft Whitney stratification\textquotedblright). (The
Whitney property is also called \emph{normal regularity} of the
stratification, see \cite[Definition 5]{Ioffe08}.) The
\emph{dimension }$\dim\left(  X\right)  $ of a semialgebraic set $X$
can thus be defined as the dimension of the manifold of highest
dimension of its stratification. This dimension is well defined and
independent of the stratification of $X$
\cite[Section~3.3]{Coste99}.

\bigskip

As a matter of the fact, semialgebraic sets constitute an
\emph{o-minimal structure}. Let us recall the definitions of the
latter\ (see for instance \cite{Coste02},
\cite{vdDries98}).\smallskip

\begin{definition}
[o-minimal structure]\label{Definition_o-minimal} An o-minimal
structure on $(\mathbb{R},+,.)$ is a sequence of Boolean algebras
$\mathcal{O} =\{\mathcal{O}_{n}\}$, where each algebra
$\mathcal{O}_{n}$ consists of subsets of $\mathbb{R}^{n}$, called
\emph{definable} (in $\mathcal{O}$), and such that for every
dimension $n\in\mathbb{N}$ the following properties hold.

\begin{enumerate}
\item[(i)] For any set $A$ belonging to $\mathcal{O}_{n}$, both $A\times
\mathbb{R}$ and $\mathbb{R}\times A$ belong to $\mathcal{O}_{n+1}$.

\item[(ii)] If $\Pi:\mathbb{R}^{n+1}\rightarrow\mathbb{R}^{n}$ denotes the
canonical projection, then for any set $A$ belonging to $\mathcal{O}_{n+1}$,
the set $\Pi(A)$ belongs to $\mathcal{O}_{n}$.

\item[(iii)] $\mathcal{O}_{n}$ contains every set of the form $\{x\in
\mathbb{R}^{n} : p(x)=0\}$, for polynomials $p:\mathbb{R}^{n}\rightarrow
\mathbb{R}$.

\item[(iv)] The elements of $\mathcal{O}_{1}$ are exactly the finite unions of
intervals and points.
\end{enumerate}

When $\mathcal{O}$ is a given o-minimal structure, a function
$f:\mathbb{R} ^{n}\rightarrow\mathbb{R}^{m}$ (or a set-valued
mapping $F:\mathbb{R} ^{n}\rightrightarrows\mathbb{R}^{m}$) is
called \emph{definable} (in $\mathcal{O}$) if its graph is definable
as a subset of $\mathbb{R}^{n} \times\mathbb{R}^{m}$.
\end{definition}

It is obvious by definition that semialgebraic sets are stable under
Boolean operations. As a consequence of the Tarski-Seidenberg
principle, they are also stable under projections, thus they satisfy
the above properties. Nonetheless, broader o-minimal structures also
exist. In particular, the Gabrielov theorem implies that
\textquotedblleft globally subanalytic\textquotedblright\ sets are
o-minimal. These two structures in particular provide rich practical
tools, because checking semi-algebraicity or subanalyticity of sets
in concrete problems of variational analysis is often easy. We refer
to \cite{BDL04}, \cite{BDL2008}, and \cite{Ioffe2009} for details.
Let us mention that Theorem~\ref{theorem:stratification} still holds
in an arbitrary o-minimal structure (it is sufficient to replace the
word \textquotedblleft semialgebraic\textquotedblright\ by
\textquotedblleft definable\textquotedblright\ in the statement). As
a matter of the fact, the statement of
Theorem~\ref{theorem:stratification} can be reinforced for definable
sets (namely, the stratification can be taken analytic), but this is
not necessary for our purposes.

\bigskip

\noindent\textbf{Remark.} Besides formulating our results and main
theorem for semialgebraic sets (the reason being their simple
definition), the validity of these results is not confined to this
class. In fact, all forthcoming statements will still hold for
\textquotedblleft definable\textquotedblright \ sets (replace
\textquotedblleft semialgebraic\textquotedblright\ by
\textquotedblleft definable in an o-minimal
structure\textquotedblright) with an identical proof. Moreover,
since our key arguments are essentially of a local nature, one can
go even further and formulate the results for the so-called
\emph{tame} sets (e.g. \cite{BDL2008}, \cite{Ioffe2009}), that is,
sets whose intersection with every ball is definable in some
o-minimal structure. (In the latter case though, slight technical
details should be taken into consideration.)

\medskip

We close this section by mentioning an important property of
semialgebraic (more generally, o-minimal) sets. Let us recall that
(topological) genericity and full measure (\textit{i.e.,} almost
everywhere) are different ways to affirm that a given property holds
in a large set. However, these notions are often complementary. In
particular, it is possible for a (topologically) generic subset of
$\mathbb{R}^{n}$ to be of null measure, or for a full measure set to
be meager (see \cite{Oxtoby} for example). Nonetheless, this
situation disappears in our setting.

\begin{proposition}
[Genericity in a semialgebraic setting]\label{Prop:generic} Let $U,V$ be
semialgebraic subsets of $\mathbb{R}^{n},$ and assume $V\subset U.$ Then the
following properties are equivalent:\smallskip

(i) $V$ is dense in $U$ ;

(ii) $V$ is (topologically) generic in $U$ ;

(iii) $U\setminus V$ is of null (Lebesgue) measure ;

(iv) the dimension of $U\setminus V$ is strictly smaller than that
of $U$.
\end{proposition}

\subsection{Some technical results}

\label{subsec-2}

In the sequel we shall always consider a set-valued map
$S:\mathcal{X} \rightrightarrows\mathbb{R}^{m}$, where
$\mathcal{X}\subset\mathbb{R}^{n},$ and we shall assume that $S$ is
semialgebraic.

\begin{theorem}
\label{theorem:generic-strict-cty}Assume that
$S:\mathcal{X}\rightrightarrows \mathbb{R}^{m}$ is outer
semicontinuous, and the sets $\mathcal{X} \subset\mathbb{R}^{n}$ and
$\mbox{\rm Graph}\left(  S\right)  $ are semi-algebraic. Then $S$ is
strictly continuous with respect to $\mathcal{X}$ everywhere except
on a set of dimension at most $\left(  \dim\mathcal{X} -1\right)  $.
\end{theorem}

\noindent\textit{Proof.} Using Theorem~\ref{theorem:stratification}
we stratify $\mathcal{X}$ into a disjoint union of manifolds
(strata) $\{\mathcal{X}_{j}\}_{j}$ and study how $S$ behaves on the
strata $\mathcal{X}_{j}$ of full dimension (that is, $\dim\left(
\mathcal{X} _{j}\right)  =\dim\left(  \mathcal{X}\right) =d\le n $).
For each such stratum $\mathcal{X}_{j}$, if $S$ is not strictly
continuous at $\bar{x}\in \mathcal{X}_{j}$ relative to
$\mathcal{X}_{j}$, then by \cite[Theorem~9.38]{RW98}, there is some
$\bar{y}\in S(\bar{x})$ such that $\mbox{\rm
lip}_{\mathcal{X}_{j}}S\left(\bar{x}\mid\bar{y}\right) =\infty$.
Since $S$ is outer semicontinuous, we deduce from
Proposition~\ref{pro:Mordukhovich-extended} that there is a nonzero
vector $v\in N_{{\scriptsize \mbox{\rm Graph}\left(  S\right)
}}\left(  \bar {x},\bar{y}\right)  \cap\left(
T_{\mathcal{X}_{j}}\left(  \bar{x}\right) \times\left\{
\mathbf{0}_{m}\right\}  \right)$. \smallskip\newline

We now stratify the semialgebraic set $\mbox{\rm
Graph}\left(S\right) \cap\left(\mathcal{X}_{j}
\times\mathbb{R}^{m}\right)$ into a finite union of disjoint
manifolds $\{\mathcal{M}_{k}\}_{k}$. Since $v\in N_{{\scriptsize
\mbox{\rm Graph}\left( S\right)}}\left( \bar{x},\bar {y}\right)
\setminus \{\mathbf{0}_{n+m}\}$, it can be written as a limit of
Hadamard normal vectors $v_{i}\in\hat{N}_{{\scriptsize \mbox{\rm
Graph}\left( S\right) }}\left( x_{i},y_{i}\right)  $ with $\left(
x_{i},y_{i}\right) \rightarrow\left( \bar{x},\bar{y}\right) $.
Passing to a subsequence if necessary, we may assume that the
sequence $\{(x_{i},y_{i})\}_{i}$ belongs to the same stratum, say
$\mathcal{M}_{k^*}\ $ and $v_{i}\in\hat{N}_{{\scriptsize \mathcal{M}
_{k^*}}}\left( x_{i},y_{i}\right) $ (note that $\mathcal{M}_{k^*}
\subset\mbox{\rm Graph}\left(  S\right)$). Since $\mathcal{M}_{k^*}$
is a smooth manifold, we have $\hat{N}_{{\scriptsize
\mathcal{M}_{k^*}}}\left( x_{i},y_{i}\right) =N_{{\scriptsize
\mathcal{M}_{k^*}}}\left( x_{i} ,y_{i}\right) =[T_{{\scriptsize
\mathcal{M}_{k^*}}}\left( x_{i} ,y_{i}\right) ]^{\bot}$. Using the
Whitney property (normal regularity) of the stratification, we
deduce that $v$ must lie in some $N_{{\scriptsize
\mathcal{M}}}\left(  \bar{x},\bar{y}\right) \cap\left(
T_{\mathcal{X}_{j}}\left(  \bar{x}\right)  \times\left\{
\mathbf{0}_{m}\right\} \right)$, where $\mathcal{M}$ is the stratum
that contains $\left(  \bar {x},\bar{y}\right)$.
Lemma~\ref{lemma:Ioffes-lemma} then tells us that the set of all
possible $\bar{x}$ is of lower dimension than that of
$\mathcal{X}_{j}$ (or $\mathcal{X}$). Since there are finitely many
strata $\mathcal{X}_{j}$, the result follows. $\hfill\Box$

\bigskip

\noindent\textbf{Remark.} Note that the domain of $S$
\[
\mathrm{dom}(S):=\{x\in\mathcal{X}:S(x)\neq\emptyset\},
\]
being the projection to $\mathbb{R}^{n}$ of the semialgebraic set
$\mbox{\rm Graph}\left(  S\right)  $, is always semialgebraic. Thus, if $S$
has nonempty values, the above assumption \textquotedblleft$\mathcal{X}$
semialgebraic\textquotedblright\ becomes superfluous. In any case, one can
eliminate this assumption from the statement and replace $\mathcal{X}$ by
$\mathcal{X}^{\prime}:=\mathrm{dom}(S)$ the domain of $S.$

$\bigskip$

The next lemma will be crucial in the sequel. We shall first need
some notation. In the sequel we denote by
\begin{equation}
\mathcal{L}:=\left\{  \mathbf{0}_n\right\}  \times\mathbb{R}^{m}
\label{eq:emph-L}
\end{equation}
as a subspace of $\mathbb{R}^{n}\times\mathbb{R}^{m}$ and we denote
by $\bar{S}:\mathbb{R}^{n}\rightrightarrows\mathbb{R}^{m}$ the
set-valued map whose graph is the closure of the graph of $S,$ that
is,
\[
\mbox{\rm Graph}\left(  \bar{S}\right)  =\mbox{\rm cl}\left(
\mbox{\rm Graph}\left(  S\right)  \right)  .
\]

\begin{lemma}
\label{lemma:decomposition}Let $S:\mathbb{R}^{n}\rightrightarrows
\mathbb{R}^{m}$ be a closed-valued semialgebraic set-valued map. For
any $k>0$, there is a $\mathcal{C}^{k}$ stratification
$\{\mathcal{M}_{i}\}_{i}$ of $\mbox{\rm Graph}\left(  S\right)  $
such that if $S\left(  \bar{x}\right) \neq\bar{S}\left(
\bar{x}\right)  $ for some $\bar{x}\in\mathbb{R}^{n}$, then there
exist $\bar{y}\in\mathbb{R}^{m}$, a stratum $\mathcal{M}_{i}$ of the
stratification of $\mbox{\rm Graph}\left(  S\right)  $ and a
neighborhood $U$ of $\left(  \bar{x},\bar{y}\right)  $ such that
$\left(  \bar{x},\bar {y}\right)  \in\mbox{\rm cl}\left(
\mathcal{M}_{i}\right)  $ and $$\left( \left( \bar{x},\bar{y}\right)
+\mathcal{L}\right)  \cap\mathcal{M}_{i}\cap U=\emptyset.$$
\end{lemma}

\noindent\textit{Proof.} By Theorem~\ref{theorem:stratification} we
stratify $\mbox{\rm Graph}\left(  S\right)  $ into a disjoint union
of finitely many manifolds, that is $\mbox{\rm Graph}\left( S\right)
=\cup_{i}\mathcal{M} _{i}$. Consider the set-valued map
$S_{i}:\mathbb{R}^{n}\rightrightarrows \mathbb{R}^{m}$ whose graph
consists of the manifold $\mathcal{M}_{i}$. Let further
$\dot{S}_{i}:\mathbb{R}^{n}\rightrightarrows\mathbb{R}^{m}$ be the
map such that $\dot{S}_{i}\left(  x\right)  =\mbox{cl}\left(
S_{i}\left( x\right)  \right)  $ for all $x$, and
$\bar{S}_{i}:\mathbb{R}^{n} \rightrightarrows\mathbb{R}^{m}$ be the
map whose graph is $\mbox{cl }\left( \mbox{\rm Graph}\left(
S_{i}\right)  \right)  $, also equal to $\mbox{\rm cl}\left(
\mbox{\rm Graph}\left(  \dot{S}_{i}\right)  \right)  $. Both
$\dot{S}_{i}$ and $\bar{S}_{i}$ are semialgebraic (for example,
\cite{Coste99}), and there exists a stratification of $\mbox{\rm
cl}\left( \mbox{\rm Graph}\left(  S\right)  \right)  $ such that the
graphs of $S_{i}$, $\dot{S}_{i}$ and $\bar{S}_{i}$ can be
represented as a finite union of strata of that stratification, by
Theorem~\ref{theorem:stratification} again.\smallskip\newline We now
prove that if $S\left(  \bar{x}\right) \neq\bar{S}\left(
\bar{x}\right)  $, then there is some $i$ such that $\dot{S}_{i}$ is
not outer semicontinuous at $\bar{x}$. Indeed, in this case there
exists $\bar{y}$ such that $\left(  \bar{x},\bar{y}\right)
\in\mbox{\rm cl}\left(  \mbox{\rm Graph}\left(  S\right)  \right)
\backslash\mbox{\rm Graph}\left(  S\right)$. Note that $\mbox{\rm
cl}\left( \mbox{\rm Graph}\left(  S\right)  \right)
=\cup_{i}\mbox{\rm Graph}\left( \bar{S}_{i}\right)  $. This means
that $\left(  \bar{x},\bar{y}\right)  $ must lie in $\mbox{\rm
Graph}\left(  \bar{S}_{i}\right)  \backslash \mbox{\rm Graph}\left(
\dot{S}_{i}\right)  $ for some $i$, which means that $\dot{S}_{i}$
is not outer semicontinuous at $\bar{x}$ as claimed.\smallskip
\newline Obviously $\left(  \bar{x},\bar{y}\right)  \in\mbox{\rm cl}\left(
\mathcal{M}_{i}\right)  $. Suppose that $\left(  \left(
\bar{x},\bar {y}\right)  +\mathcal{L}\right)
\cap\mathcal{M}_{i}\cap U\neq\emptyset$ for all neighborhoods $U$
containing $\left(  \bar{x},\bar{y}\right)  $. Then there is a
sequence $y_{j}\rightarrow\bar{y}$ such that $\left(  \bar{x}
,y_{j}\right)  \in\mathcal{M}_{i}$. Since $\dot{S}_{i}$ is
closed-valued, this would yield $\left(  \bar{x},\bar{y}\right)
\in\mbox{\rm Graph}\left( \dot{S}_{i}\right)  $, which contradicts
$\left(  \bar{x},\bar{y}\right) \notin\mbox{\rm Graph}\left(
\dot{S}_{i}\right)  $ earlier. $\hfill\Box$

\bigskip

Keeping now the notation of the proof of the previous lemma, let us
set $\bar{z}:=\left(  \bar{x},\bar{y}\right)  $. Let further
$\mathcal{M} _{i},\mathcal{M}^{\prime}$ be the strata of $\mbox{\rm
cl}\left( \mbox{\rm Graph}\left(  S\right)  \right)  $ such that
$z\in\mathcal{M} ^{\prime}\subset\mbox{cl}\left(
\mathcal{M}_{i}\right)  $. In the next lemma we are working with
normals on manifolds, so it does not matter which kind of normal in
Definition~\ref{def:Normal-cones} we consider.

\begin{lemma}
\label{lemma:normals-behave}Suppose there is a neighborhood $U$ of
$\bar{z}$ such that $\bar{z}\in\mathcal{M}^{\prime}$,
$\mathcal{M}^{\prime} \subset\mbox{\rm cl}\left(
\mathcal{M}_{i}\right)  $ and $\left(  \bar {z}+\mathcal{L}\right)
\cap\mathcal{M}_{i}\cap U=\emptyset$, where $\mathcal{L}$ is defined
in \eqref{eq:emph-L}. Then $N_{\mathcal{M}^{\prime} }\left(
\bar{z}\right)  \cap\mathcal{L}^{\perp}\supsetneq\left\{
\mathbf{0}_{n+m}\right\}  $.
\end{lemma}

\noindent\textit{Proof.} We prove the result by contradiction.
Suppose that $N_{\mathcal{M}^{\prime}}\left( \bar{z}\right)
\cap\mathcal{L}^{\perp }=\left\{  \mathbf{0}_{n+m}\right\}$. Then
$T_{\mathcal{M}^{\prime}}\left( \bar{z}\right)
+\mathcal{L}=\mathbb{R}^{n+m}$ by
Proposition~\ref{pro:orthogonal-intersection-sum}. We may assume, by
taking a submanifold of $\mathcal{M}^{\prime}$ if necessary, that
$\dim\mathcal{M} ^{\prime}=n$ so that
$\dim\mathcal{M}^{\prime}+\dim\mathcal{L}=n+m$ and
$T_{\mathcal{M}^{\prime}}\left(  \bar{z}\right)  \oplus\mathcal{L}
=\mathbb{R}^{n+m}$. Owing to the so-called wink lemma (see
\cite[Proposition~5.10]{Denkowska} {\it e.g.}) we may assume that
$\dim\mathcal{M}_{i}=n+1$.\medskip

\textbf{(Case }$m=1$\textbf{)} We first consider the case where
$m=1$. {In this case,} the subspace $\mathcal{L}$ is a line whose spanning
vector $v=\left( \mathbf{0},1\right)  $ is not in
$T_{\mathcal{M}^{\prime}}\left( \bar {z}\right)  $. There is a
neighborhood $U^{\prime}$ of $\bar{z}$ such that $U^{\prime}\subset
U$, $\mathcal{M}^{\prime}\cap U^{\prime}$ equals $f^{-1}\left(
0\right)  $ for some smooth function $f:U^{\prime}
\rightarrow\mathbb{R}$ (local equation of $\mathcal{M}^{\prime}$),
and $\mathcal{M}_{i}\cap U^{\prime}=f^{-1}\left(  \left(
0,\infty\right) \right)  $. The gradient $\nabla f\left(
\bar{z}\right)  $ is nonzero and is not orthogonal to $v$ since
$T_{\mathcal{M}^{\prime}}\left( \bar{z}\right)  $ is the set of
vectors orthogonal to $\nabla f\left(  \bar{z}\right)  $ and
$T_{\mathcal{M}^{\prime}}\left( \bar{z}\right)  \oplus\mathcal{L}
=\mathbb{R}^{n+1}$. There are points in $\left(
\bar{z}+\mathcal{L}\right) \cap U^{\prime}$ such that $f$ is
positive, which means that $\left(  \bar {z}+\mathcal{L}\right)
\cap\mathcal{M}_{i}\cap U^{\prime}\neq\emptyset$, contradicting the
stipulated conditions. Therefore, we assume that $m>1$ for the rest
of the proof.

\smallskip

\textbf{(Case} $m>1$\textbf{)} As in the previous case, we shall
eventually prove that
$\left(\bar{z}+\mathcal{L}\right)\cap\mathcal{M}_{i}\cap
U^{\prime}\neq\emptyset$ reaching to a contradiction. To this end,
let us denote by $h_0$ the (semialgebraic) homeomorphism of
$\mathbb{R}^{n+m}$ to $\mathbb{R}^{n+m}$ which, for some
neighborhood $V\subset U$ of $\bar{z}$, maps homeomorphically
$V\cap\left(\mathcal{M}_{i}\cup\mathcal{M}^{\prime}\right)$ to
$\mathbb{R}^{n}\times\left(\mathbb{R}_{+}\times\left\{
\mathbf{0}_{m-1}\right\}\right) \subset\mathbb{R}^{n+m}$ and $V \cap
\mathcal{M}^\prime$ to $\mathbb{R}^n\times\{\mathbf{0}_m\}$ (see
\cite[Theorem~3.12]{Coste99} {\it e.g.}). \smallskip

{\it Claim.} We first show that there exists a closed neighborhood
$W\subset V$ of $\bar{z}$ such that $W\cap\mathcal{M}^{\prime}$ and
$\partial W\cap\mathcal{M}_{i}$ are both homeomorphic to
$\mathbb{B}^{n}$ and
$W\cap\mathcal{M}^{\prime}=\mathbb{B}^{n+m}\left(\bar{z},R_1\right)\cap\mathcal{M}^{\prime}$
for some $R_{1}>0$.\smallskip

\noindent {Since} $\mathcal{M}^{\prime}$ is a smooth manifold,
there exists $R_{1}>0$ such that
$\mathbb{B}^{n+m}\left(\bar{z},R_{1}\right)\cap\mathcal{M}^{\prime}$
is homeomorphic (in fact, diffeomorphic) to
$(T_{\mathcal{M}^\prime}(\bar{z})+\bar{z})\cap
\mathbb{B}^{n+m}(\bar{z},R_1)$, which in turn is homeomorphic to
$\mathbb{B}^{n}$, as is shown by the homeomorphism:
\[
z\mapsto\left(\frac{\left\vert z-\bar{z}\right\vert }{\left\vert
P\left(z\right)-\bar{z}\right\vert
}\left(P\left(z\right)-\bar{z}\right)\right)+\bar{z}\,,
\]
where $P$
denotes the projection onto the tangent space
$\bar{z}+T_{\mathcal{M}^{\prime}}\left(\bar{z}\right)$. Consider the
image of
$\mathbb{B}^{n+m}\left(\bar{z},R_{1}\right)\cap\mathcal{M}^{\prime}$
under the map $h_{0}$. This image lies in the set
$\mathbb{R}^{n}\times\left\{ \mathbf{0}_m\right\}$. Therefore, for
$r_{1}>0$ sufficiently small, the set $W=
h_{0}^{-1}\left(h_{0}\left(\mathbb{B}^{n+m}\left(\bar{z},R_{1}\right)\cap\mathcal{M}^{\prime}\right)+\left[-r_{1},r_{1}\right]^{m}\right)$
satisfies the required properties, {concluding the proof of our claim}. \smallskip

\begin{figure}
\begin{longtable}{ccc}
\begin{tabular}{c}
\includegraphics[scale=0.3]{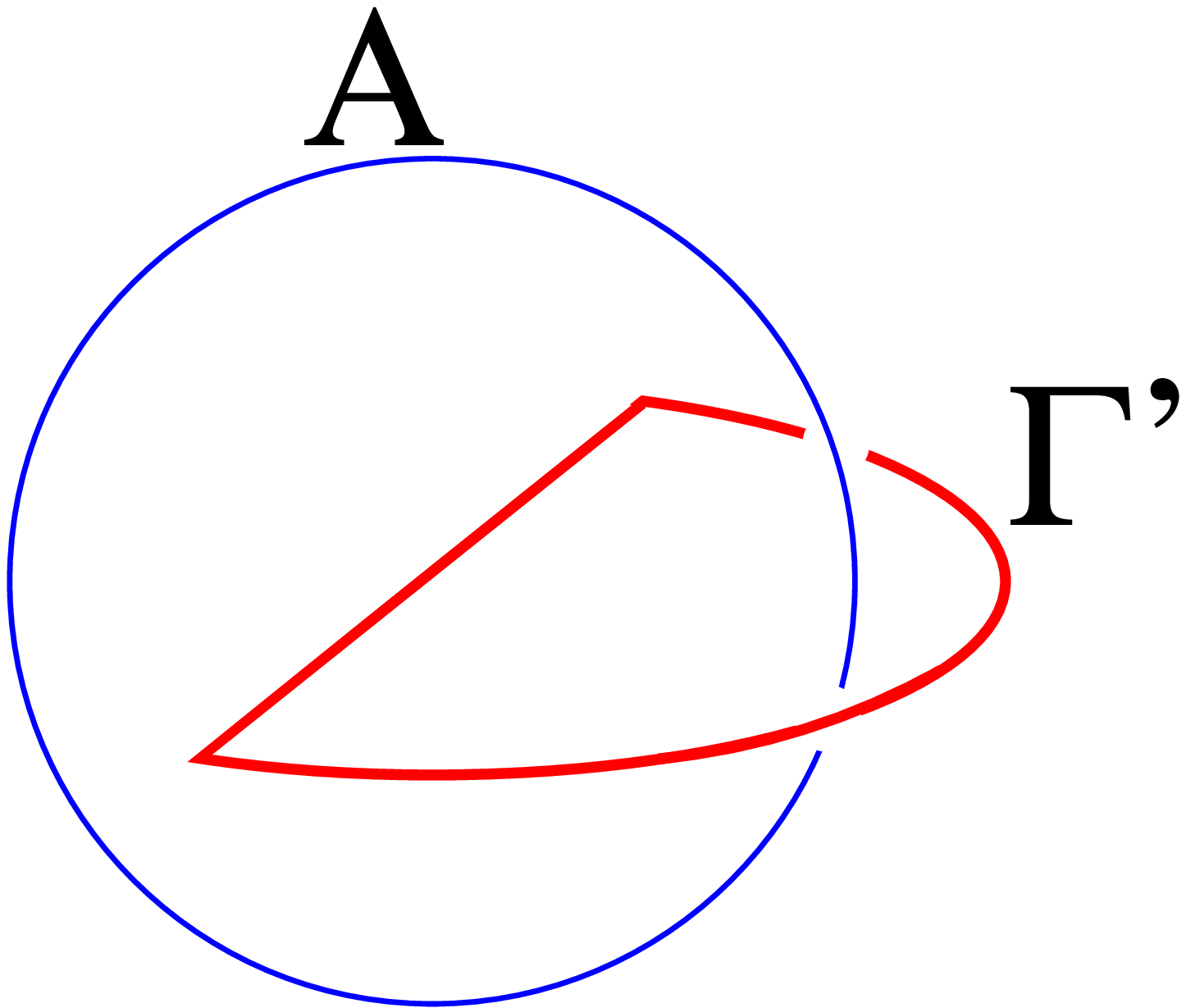}\tabularnewline
\end{tabular} & {\Huge $\rightarrow$}  & \begin{tabular}{c}
\includegraphics[scale=0.3]{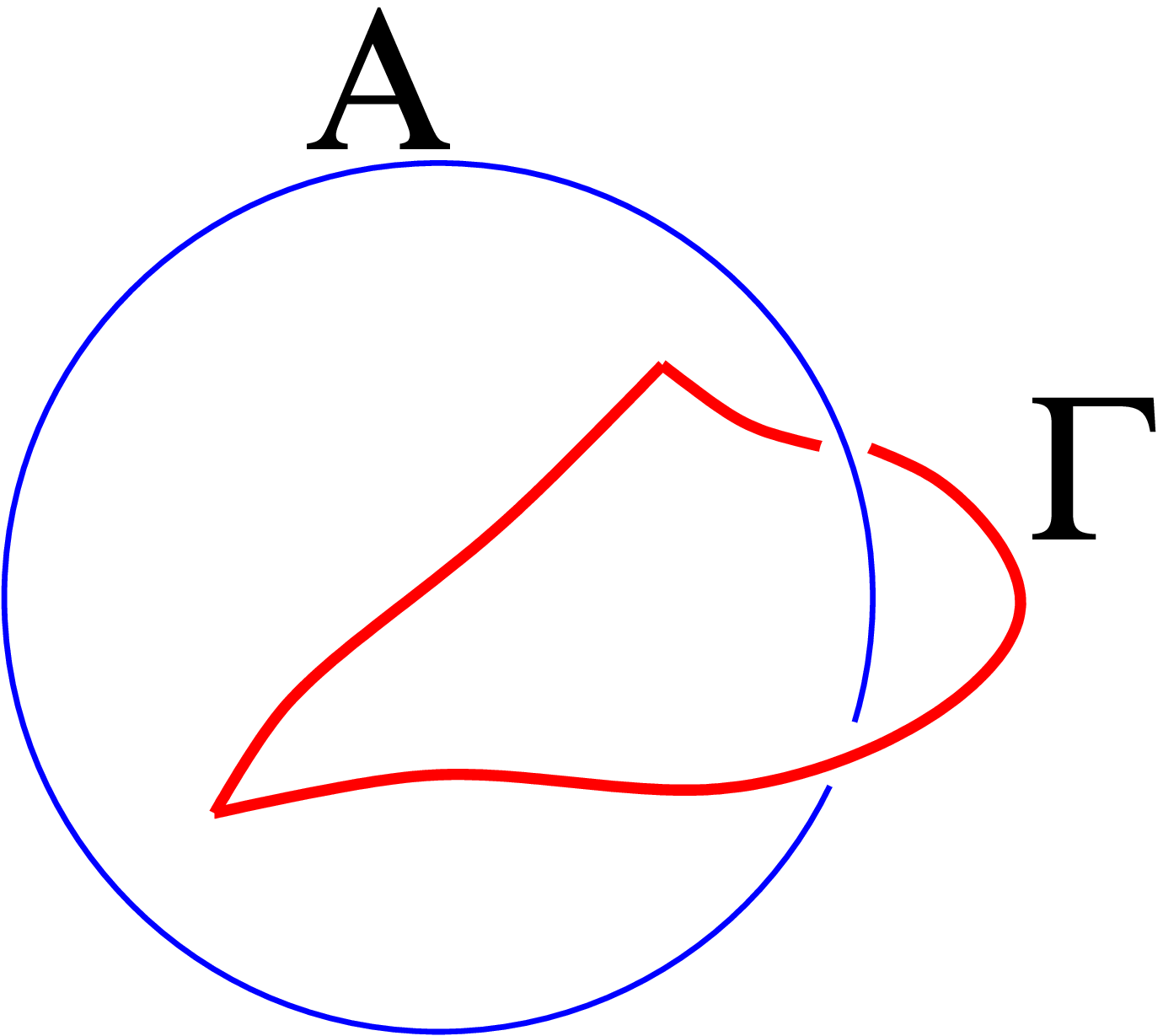}\tabularnewline
\end{tabular}\tabularnewline
\end{longtable}
\caption{\label{fig:Link-fig} Linking sets
$\left(A,\Gamma^{\prime}\right)$ and $\left(A,\Gamma\right)$.}
\end{figure}

Let us further fix $v\in\mathcal{L}\backslash\left\{
\mathbf{0}\right\}$ and consider the set
\[
\Gamma^{\prime}:=\underbrace{\left(\mathbb{B}^{n+m}\left(\bar{z},R\right)\cap\left(\bar{z}+T_{\mathcal{M}^{\prime}}\left(\bar{z}\right)\right)\right)}_{\Gamma_{1}^{\prime}}\cup\underbrace{\left(\mathbb{S}^{n+m-1}\left(\bar{z},R\right)\cap\left(\bar{z}+T_{\mathcal{M}^{\prime}}\left(\bar{z}\right)+\mathbb{R}_{+}\left\{
v\right\} \right)\right)}_{\Gamma_{2}^{\prime}}\,.
\]
Setting
\[
A:=\mathbb{S}^{n+m-1}\left(\bar{z},r\right)\cap\mathcal{L},\mbox{
where }0<r<R\,,
\] we immediately get that the sets $A$ and
$\Gamma^{\prime}$ link ({\it c.f.}
Theorem~\ref{theorem:linking-sets}). Based on this, our objective is
to prove that the sets $A$ {and $\Gamma$ also link, where $\Gamma$
is} defined by
\[
\Gamma=\underbrace{\left(W\cap\mathcal{M}^{\prime}\right)}_{\Gamma_{1}}\cup\underbrace{\left(\partial
W\cap(\mathcal{M}_{i}\cup
\mathcal{M}^\prime)\right)}_{\Gamma_{2}}\,,
\]
provided $r>0$ is chosen appropriately. Once we succeed in doing so,
we apply Definition~\ref{Def:link} (for $h=id$) to deduce that
$\left(\bar{z}+\mathcal{L}\right)\cap\mathcal{M}_{i}\cap
U\neq\emptyset$, which contradicts our initial assumptions.
\newline Figure~\ref{fig:Link-fig} illustrates the sets $A$, $\Gamma$ and
$\Gamma^{\prime}$ for $n=1$ and $m=2$. \smallskip

For the sequel, we introduce the notation
{}``$\xrightarrow{\simeq}{}$\textquotedblright\ in
$f:D_{1}\xrightarrow{\simeq}{}D_{2}$ to mean that $f$ is a
homeomorphism between the sets $D_{1}$ and $D_{2}$. In Step~1 and
Step~2, we define a continuous function
$H:\left(\mathbb{B}^{n+1}\left(\mathbf{0},1\right)\times\left\{
0\right\}
\right)\cup\left(\mathbb{S}^{n}\left(\mathbf{0},1\right)\times\left[0,2\right]\right)\rightarrow\mathbb{B}^{n+m}\left(\bar{z},R\right)$
that will be used in Step~3. \medskip

\textbf{Step~1: Determine $H$ on
$\left(\mathbb{B}^{n+1}(\mathbf{0},1)\times\left\{ 0\right\}
\right)\cup\left(\mathbb{S}^{n}\left(\mathbf{0},1\right)\times\left[0,2\right]\right)$}.

In Steps~1\,($a$) to 1\,($c$), we define a continuous function $H$
on $\mathbb{S}^{n}\left(\mathbf{0},1\right)\times\left[0,2\right]$
so that
$H\mid_{\mathbb{S}^{n}\left(\mathbf{0},1\right)\times\left[0,2\right]}$
is a homotopy between $\Gamma$ and $\Gamma^{\prime}$. More
precisely, denoting by \begin{align*}
\mathbb{S}_{+}^{n}\left(\mathbf{0},1\right) & :=\mathbb{S}^{n}\left(\mathbf{0},1\right)\cap\left(\mathbb{R}^{n}\times[0,\infty)\right),\\
\mathbb{S}_{-}^{n}\left(\mathbf{0},1\right) & :=\mathbb{S}^{n}\left(\mathbf{0},1\right)\cap\left(\mathbb{R}^{n}\times(-\infty,0]\right),\end{align*}
 we want to define $H$ in such a way that its restrictions
 \begin{align*}
H\left(\cdot,0\right)\mid_{\mathbb{S}_{+}^{n}\left(\mathbf{0},1\right)}:\mathbb{S}_{+}^{n}\left(\mathbf{0},1\right) & \xrightarrow{\simeq}\Gamma_{1}\subset\mathbb{R}^{n+m},\\
H\left(\cdot,0\right)\mid_{\mathbb{S}_{-}^{n}\left(\mathbf{0},1\right)}:\mathbb{S}_{-}^{n}\left(\mathbf{0},1\right) & \xrightarrow{\simeq}\Gamma_{2}\subset\mathbb{R}^{n+m},\\
H\left(\cdot,2\right)\mid_{\mathbb{S}_{+}^{n}\left(\mathbf{0},1\right)}:\mathbb{S}_{+}^{n}\left(\mathbf{0},1\right) & \xrightarrow{\simeq}\Gamma_{1}^{\prime}\subset\mathbb{R}^{n+m},\\
H\left(\cdot,2\right)\mid_{\mathbb{S}_{-}^{n}\left(\mathbf{0},1\right)}:\mathbb{S}_{-}^{n}\left(\mathbf{0},1\right) & \xrightarrow{\simeq}\Gamma_{2}^{\prime}\subset\mathbb{R}^{n+m},
\end{align*}
 are homeomorphisms between the respective spaces. Note that both
$\mathbb{S}_{+}^{n}\left(\mathbf{0},1\right)$ and
$\mathbb{S}_{-}^{n}\left(\mathbf{0},1\right)$ are homeomorphic to
$\mathbb{B}^{n}\left(\mathbf{0},1\right)$. For notational
convenience, we denote by
$\mathbb{S}_{=}^{n}\left(\mathbf{0},1\right)$ the set
$\mathbb{S}^{n}\left(\mathbf{0},1\right)\cap\left(\mathbb{R}^{n}\times\left\{
0\right\}
\right)=\mathbb{S}^{n-1}\left(\mathbf{0},1\right)\times\left\{
0\right\}$.
\smallskip

\textit{Step 1\,($a$)}.\, \textbf{Determine $H$ on
$\mathbb{S}\left(\mathbf{0},1\right)\times\left[0,1\right]$}.
\\ Since $\partial W\cap \mathrm{cl}\,\mathcal{M}_{i}$ is a closed set that does
not contain $\bar{z}$, there is some $R>0$ such that $\left(\partial
W\cap\mathcal{M}_{i}\right)\cap\mathbb{B}^{n+m}\left(\bar{z},R\right)=\emptyset$
and $\mathbb{B}^{n+m}\left(\bar{z},R\right)\subset U$. We proceed to
create the homotopy $H$ so that\begin{align*}
H\left(\cdot,1\right)\mid_{\mathbb{S}_{+}^{n}\left(\mathbf{0},1\right)}:\mathbb{S}_{+}^{n}\left(\mathbf{0},1\right) & \xrightarrow{\simeq}\Gamma_{1}^{\prime\prime}\subset\mathbb{R}^{n+m},\\
H\left(\cdot,1\right)\mid_{\mathbb{S}_{-}^{n}\left(\mathbf{0},1\right)}:\mathbb{S}_{-}^{n}\left(\mathbf{0},1\right) & \xrightarrow{\simeq}\Gamma_{2}^{\prime\prime}\subset\mathbb{R}^{n+m},\end{align*}
where\begin{align*}
\Gamma_{1}^{\prime\prime} & =\mathbb{B}^{n+m}\left(\bar{z},R\right)\cap\mathcal{M}^{\prime},\\
\mbox{ and }\qquad \Gamma_{2}^{\prime\prime} &
\subset\mathbb{S}^{n+m-1}\left(\bar{z},R\right)\mbox{ is
homeomorphic to }\Gamma_{2}.\end{align*} The first homotopy between
$\Gamma_{1}$ and $\Gamma_{1}^{\prime\prime}$ can be chosen such that
$H\left(s,t\right)\in\mathcal{M}^{\prime}$ for all
$s\in\mathbb{S}_{+}^{n}\left(\mathbf{0},1\right)$ and
$t\in\left[0,1\right]$. We also require that
$d\left(\bar{z},H\left(s,t\right)\right)\geq R$ for all
$s\in\mathbb{S}_{=}^{n}\left(\mathbf{0},1\right)$ and
$t\in\left[0,1\right]$, which does not present any difficulties.

For the second homotopy between $\Gamma_{2}$ and
$\Gamma_{2}^{\prime\prime}$, we first extend $H\left(\cdot,1\right)$
so that
$H\left(\cdot,1\right)\mid_{\mathbb{S}^{n}\left(\mathbf{0},1\right)}:
\mathbb{S}^{n}\left(\mathbf{0},1\right)\xrightarrow{\simeq}{}\Gamma_{1}^{\prime\prime}\cup\Gamma_{2}^{\prime\prime}$
is a homeomorphism between the corresponding spaces. This is
achieved by showing that there is a homeomorphism
$H\left(\cdot,1\right)\mid_{\mathbb{S}_{-}^{n}\left(\mathbf{0},1\right)}$
between $\mathbb{S}_{-}^{n}\left(\mathbf{0},1\right)$ and
$\Gamma_{2}^{\prime\prime}$. Let
$h_{2}:\mathbb{B}^{n}\left(\mathbf{0},1\right)\xrightarrow{\simeq}{}\mathbb{S}_{-}^{n}\left(\mathbf{0},1\right)$
be a homeomorphism between $\mathbb{B}^{n}\left(\mathbf{0},1\right)$
and $\mathbb{S}_{-}^{n}\left(\mathbf{0},1\right)$. Then
$H\left(\cdot,1\right)\mid_{\mathbb{S}_{=}^{n}\left(\mathbf{0},1\right)}\circ
h_{2}\mid_{\mathbb{S}^{n-1}\left(\mathbf{0},1\right)}:\mathbb{S}^{n-1}\left(\mathbf{0},1\right)\xrightarrow{\simeq}{}\partial\Gamma_{2}^{\prime\prime}$.
By Lemma~\ref{lemma:extension-to-interior} this can be extended to a
homeomorphism
$G:\mathbb{B}^{n}\left(\mathbf{0},1\right)\xrightarrow{\simeq}{}\Gamma_{2}^{\prime\prime}$.
Define
$H\left(\cdot,1\right)\mid_{\mathbb{S}_{-}^{n}\left(\mathbf{0},1\right)}:\mathbb{S}_{-}^{n}\left(\mathbf{0},1\right)\xrightarrow{\simeq}{}\Gamma_{2}^{\prime\prime}$
by
$H\left(\cdot,1\right)\mid_{\mathbb{S}_{-}^{n}\left(\mathbf{0},1\right)}=G\circ
h_{2}^{-1}$.

It remains to resolve $H$ on
$\mathbb{S}_{-}^{n}\left(\mathbf{0},1\right)\times\left(0,1\right)$.
Note that the sets
$$H\left(\mathbb{S}_{=}^{n}\left(\mathbf{0},1\right)\times\left[0,1\right]\right),
\quad
H\left(\mathbb{S}_{-}^{n}\left(\mathbf{0},1\right)\times\left\{
0\right\} \right)=\Gamma_{2}\quad \text{and}\quad
H\left(\mathbb{S}_{-}^{n}\left(\mathbf{0},1\right)\times\left\{
1\right\} \right)=\Gamma_{2}^{\prime\prime}$$ are all of dimension
at most $n$, so the radial projection of these sets onto
$\mathbb{S}^{n+m-1}\left(\bar{z},R\right)$ is of dimension at most
$n$. Since $\mathbb{S}^{n+m-1}\left(\bar{z},R\right)$ is of
dimension at least $n+1$, we can find some point
$p\in\mathbb{S}^{n+m-1}\left(\bar{z},R\right)$ not lying in the
radial projections of these sets. The set
\[D:=\mathbb{R}^{n+m}\backslash(((\mathbb{R}_{+}\left\{ p-\bar{z}\right\} )+\left\{ \bar{z}\right\} )\cup\mathbb{B}^{n+m}\left(\bar{z},R\right))\]
is homeomorphic to $\mathbb{R}^{n+m}$, so by the Tietze extension theorem (see
for example \cite{Munkres00}), we can extend $H$ continuously to
$\mathbb{S}_{-}^{n}\left(\mathbf{0},1\right)\times\left[0,1\right]$
so that
$H(\mathbb{S}_{-}^{n}\left(\mathbf{0},1\right)\times\left[0,1\right])\subset
D$.

\smallskip

\textit{Step 1\,($b$)}.\,\textbf{Determine $H$ on
$\mathbb{S}_{+}^{n}\left(\mathbf{0},1\right)\times\left[1,2\right]$.}\\
We next define
$H\mid_{\mathbb{S}_{+}^{n}\left(\mathbf{0},1\right)\times\left[1,2\right]}$,
the homotopy between $\Gamma_{1}^{\prime\prime}$ and
$\Gamma_{1}^{\prime}$. Since $\mathcal{M}^{\prime}$ is a manifold,
for any $\delta>0$, we can find $R$ small enough such that for any
$z\in\mathbb{B}^{n+m}\left(\bar{z},R\right)\cap\mathcal{M}^{\prime}$,
the distance from $z$ to
$\bar{z}+T_{\mathcal{M}^{\prime}}\left(\bar{z}\right)$ is at most
$\delta R$. The value $R$ can be reduced if necessary so that the
mapping $P$, which projects a point
$z\in\mathbb{B}^{n+m}\left(\bar{z},R\right)\cap\mathcal{M}^{\prime}$
to the closest point in
$\bar{z}+T_{\mathcal{M}^{\prime}}\left(\bar{z}\right)$, is a
homeomorphism of
$\mathbb{B}^{n+m}\left(\bar{z},R\right)\cap\mathcal{M}^{\prime}$ to
its image.

Define the map
$H_{1}:\left(\mathbb{B}^{n+m}\left(\bar{z},R\right)\cap\mathcal{M}^{\prime}\right)\times\left[1,2\right]\rightarrow\mathbb{B}^{n+m}\left(\bar{z},R\right)$
by \[ H_{1}\left(z,t\right):=\left(\frac{\left\vert
z-\bar{z}\right\vert }{\left\vert
\left(2-t\right)z+\left(t-1\right)P\left(z\right)-\bar{z}\right\vert
}\left(\left(2-t\right)z+\left(t-1\right)P\left(z\right)-\bar{z}\right)\right)+\bar{z}.\]
 This is a homotopy from $\Gamma_{1}$ to $\Gamma_{1}^{\prime}$.
For any homeomorphism
$h_{1}:\mathbb{B}^{n+m}\left(\bar{z},R\right)\cap\mathcal{M}^{\prime}\xrightarrow{\simeq}{}\mathbb{S}_{+}^{n}\left(\mathbf{0},1\right)$,
we define
$H\mid_{\mathbb{S}_{+}^{n}\left(\mathbf{0},1\right)\times\left[0,1\right]}$
via
$H\left(s,t\right)=H_{1}\left(h_{1}^{-1}\left(s\right),t\right)$.

\smallskip

\textit{Step 1\,($c$)}.\,\textbf{Determine $H$ on
$\mathbb{S}_{-}^{n}\left(\mathbf{0},1\right)\times\left[1,2\right]$.}
We now define
$H\mid_{\mathbb{S}_{-}^{n}\left(\mathbf{0},1\right)\times\left[1,2\right]}$,
the homotopy between $\Gamma_{2}^{\prime\prime}$ and
$\Gamma_{2}^{\prime}$ that respects the boundary conditions
stipulated by
$H\mid_{\mathbb{S}_{=}^{n}\left(\mathbf{0},1\right)\times\left[1,2\right]}$.
We extend
$H\left(\cdot,1\right)\mid_{\mathbb{S}^{n}\left(\mathbf{0},1\right)}$
so that it is a homeomorphism between
$\mathbb{S}^{n}\left(\mathbf{0},1\right)$ and
$\Gamma_{1}^{\prime}\cup\Gamma_{2}^{\prime}$ by using methods
similar to that used in Step 1($a$).

We now use the Tietze extension theorem to establish a continuous
extension of $H$ to
$\mathbb{S}^{n}\left(\mathbf{0},1\right)\times\left[1,2\right]$. We
are left only to resolve $H$ on
$\mathbb{S}_{-}^{n}\left(\mathbf{0},1\right)\times\left(1,2\right)$.
Much of this is now similar to the end of step 1($a$). The dimension of
$\mathbb{S}^{n+m-1}\left(\bar{z},R\right)$ is $n+m-1$, while the
dimensions of $\Gamma_{2}^{\prime\prime}$, $\Gamma_{2}^{\prime}$ and
$H\left(\mathbb{S}_{=}^{n}\left(\mathbf{0},1\right)\times\left[1,2\right]\right)$
are all at most $n$. Therefore, there is one point in
$\mathbb{S}^{n+m-1}\left(\bar{z},R\right)$ outside these three sets,
say $p$. Since
$\mathbb{S}^{n+m-1}\left(\bar{z},R\right)\backslash\left\{ p\right\}
$ is homeomorphic to $\mathbb{R}^{n+m-1}$, the Tietze extension
theorem again implies that we can extend $H$ continuously in
$\mathbb{S}^{n}\left(\mathbf{0},1\right)\times\left[1,2\right]$.

\smallskip

\textit{Step 1\,($d$)}.\,\textbf{Determine $H$ on}
$\mathbb{B}^{n+1}\left(\mathbf{0},1\right)\times\left\{ 0\right\} $.
We use Lemma~\ref{lemma:extension-to-interior} to extend the domain
of the function
$$H\left(\cdot,0\right):\mathbb{S}^{n}\left(\mathbf{0},1\right)\xrightarrow{\simeq}{}
{\left(\mathcal{M}^{\prime}\cap W \right) \cup (\mathcal{M}_i \cap
\partial W)}$$
 to $\mathbb{B}^{n+1}\left(\mathbf{0},1\right)$ so that \[
H\left(\cdot,0\right):\mathbb{B}^{n+1}\left(\mathbf{0},1\right)\xrightarrow{\simeq}{}\left(\mathcal{M}^{\prime}\cup\mathcal{M}_{i}\right)\cap
W\]
 is a homeomorphism.

\smallskip

\textbf{Step 2: Choice of $R$ and $r$.} We now choose $R$ and $r$ so
that
$H\left(\mathbb{S}^{n}\left(\mathbf{0},1\right)\times\left[0,2\right]\right)$
does not intersect
$A=\mathbb{S}^{n+m-1}\left(\bar{z},r\right)\cap\left(\bar{z}+\mathcal{L}\right)$.
To this end, consider the minimization problem \[ \min\,\left\{
\mathrm{dist}\,\left(z,T_{\mathcal{M}^{\prime}}\left(\bar{z}\right)+\bar{z}\right):z\in\mathbb{S}^{n+m-1}\left(\bar{z},r\right)\cap\left(\bar{z}+\mathcal{L}\right)\right\}
.\] Since
$\mathbb{S}^{n+m-1}\left(\bar{z},r\right)\cap\left(\bar{z}+\mathcal{L}\right)$
is compact, the above minimum is attained at some point $z_{r}$ and
its value is not zero (otherwise $z_{r}-\bar{z}$ would be a nonzero
element in
$T_{\mathcal{M}^{\prime}}\left(\bar{z}\right)\cap\mathcal{L}$,
contradicting
$T_{\mathcal{M}^{\prime}}\left(\bar{z}\right)\cap\mathcal{L}=\left\{
\mathbf{0}\right\}$). Therefore, for some constant $\varepsilon\in
(0,1]$ independent of $r$, it holds
$\mathrm{dist}\,\left(z_{r},T_{\mathcal{M}^{\prime}}\left(\bar{z}\right)+\bar{z}\right)=\varepsilon\,
r$.

Given $\delta>0$, we can shrink $R$ if necessary to get
$\mbox{d}\left(z,T_{\mathcal{M}^{\prime}}\left(\bar{z}\right)+\bar{z}\right)\leq\delta\,
R$ for all $z\in
H\left(\mathbb{S}_{+}^{n}\left(\mathbf{0},1\right)\times\left[0,1\right]\right)$.
If $\delta<\varepsilon$, we can find some $r$ satisfying $\delta
R<\varepsilon\, r\leq r<R$. Since $\delta R<\varepsilon\, r$,
$H\left(\mathbb{S}_{+}^{n}\left(\mathbf{0},1\right)\times\left[1,2\right]\right)$
does not intersect
$\mathbb{S}^{n+m-1}\left(\bar{z},r\right)\cap\left(\bar{z}+\mathcal{L}\right)$.
From $r<R$, it is clear that
$H\left(\mathbb{S}_{-}^{n}\left(\mathbf{0},1\right)\times\left[0,2\right]\right)$,
being a subset of
$\mbox{cl}\left(\mathbb{R}^{n+m}\backslash\mathbb{B}^{n+m}\left(\bar{z},R\right)\right)$,
does not intersect
$\mathbb{S}^{n+m-1}\left(\bar{z},r\right)\cap\left(\bar{z}+\mathcal{L}\right)$.
Elements in
$H\left(\mathbb{S}_{+}^{n}\left(\mathbf{0},1\right)\times\left[0,1\right]\right)$
are either in $\mathbb{B}^{n+m}(\bar{z},R) \cap \mathcal{M}^\prime$
or outside $\mathbb{B}^{n+m}(\bar{z},R)$, so
$H\left(\mathbb{S}^{n}\left(\mathbf{0},1\right)\times\left[0,2\right]\right)$
does not intersect $A$ as needed.

\smallskip

\textbf{Step 3: ``Set-up\textquotedblright\ for linking theorem}.
Let
\[
h_{3}:\mathbb{B}^{n+1}\left(\mathbf{0},1\right)\xrightarrow{\simeq}{}\left(\mathbb{B}^{n+1}\left(\mathbf{0},1\right)\times\left\{
0\right\}
\right)\cup\left(\mathbb{S}^{n}\left(\mathbf{0},1\right)\times\left[0,2\right]\right)\]
 be a homeomorphism between the respective spaces. We can extend the
homeomorphism
\[
H\mid_{\mathbb{S}^{n}\left(\mathbf{0},1\right)\times\left\{
2\right\} }\circ
h_{3}\mid_{\mathbb{S}^{n}\left(\mathbf{0},1\right)}:\mathbb{S}^{n}\left(\mathbf{0},1\right)\xrightarrow{\simeq}{}\Gamma^{\prime}
\]
to
\[
h_{4}:\mathbb{B}^{n+1}\left(\mathbf{0},1\right)\xrightarrow{\simeq}{}\left(T_{\mathcal{M}^{\prime}}\left(\bar{z}\right)+\mathbb{R}_{+}\left\{
v\right\}+\bar{z}\right)\cap\mathbb{B}^{n+m}\left(\bar{z},R\right).\]
Define the map
\[
g:\left(T_{\mathcal{M}^{\prime}}\left(\bar{z}\right)+\mathbb{R}_{+}\left\{
v\right\}+\bar{z}\right)\cap\mathbb{B}^{n+m}\left(\bar{z},R\right)\rightarrow\mathbb{B}^{n+m}\left(\bar{z},R\right)
\]
by $g=H\circ h_{3}\circ h_{4}^{-1}$. By construction, the map
$g\mid_{\Gamma^{\prime}}$ is the identity map there. Furthermore,
$g$ can be extended continuously to the domain $\mathbb{R}^{n+m}$ by
the Tietze extension theorem.

\smallskip

\textbf{Step 4: Apply linking theorem}. Recall that
$A:=\mathbb{B}^{n+m}\left(\bar{z},r\right)\cap\left(\bar{z}+\mathcal{L}\right)$
and $\Gamma^{\prime}$ link by Theorem~\ref{theorem:linking-sets}.
This means that there is a nonempty intersection of
$g\left(\left(T_{\mathcal{M}^{\prime}}\left(\bar{z}\right)+\mathbb{R}_{+}\left\{
v\right\}
+\bar{z}\right)\cap\mathbb{B}^{n+m}\left(\bar{z},R\right)\right)$
with $A$. Step~2 asserts that the intersection is not in
$H\left(\mathbb{S}^{n}\left(\mathbf{0},1\right)\times\left[0,2\right]\right)$,
so the intersection lies in
$H\left(\mathbb{B}^{n+1}\left(\mathbf{0},1\right)\times\left\{
0\right\} \right)$. In other words, $A$ and $\Gamma$ link. This
means that $W\cap\mathcal{M}_{i}$ intersects
$\mathbb{B}^{n+m}\left(\bar{z},r\right)\cap\left(\bar{z}+\mathcal{L}\right)$,
which means that
$\left(\bar{z}+\mathcal{L}\right)\cap\mathcal{M}_{i}\cap
U\neq\emptyset$, contradicting our assumption. $\hfill\Box$

\subsection{Main result}

\label{subsec-3}

In this section we put together all previous results to obtain the
following theorem. Recall that $\bar{S}$ is the set-valued map whose
graph is the closure of the graph of $S$ (thus, $\bar{S}$ is outer
semicontinuous by definition).

\begin{theorem}
\label{theorem:generic-osc}If
$S:\mathcal{X}\rightrightarrows\mathbb{R}^{m}$ is a closed-valued
semialgebraic set-valued map, where $\mathcal{X}
\subset\mathbb{R}^{n}$ is semialgebraic, then $S$ and $\bar{S}$
differ outside a set of dimension at most $\left(
\dim\mathcal{X}-1\right)  $.
\end{theorem}

\noindent\textit{Proof.} We first consider the case where
$\mathcal{X} =\mathbb{R}^{n}$ and a $\mathcal{C}^{k}$ stratification
of $\mbox{\rm cl}\left(  \mbox{\rm Graph}\left(  S\right)  \right)
$. If $S\left(  \bar{x}\right)  \neq\bar{S}\left(  \bar{x}\right) $,
then Lemma~\ref{lemma:decomposition} and
Lemma~\ref{lemma:normals-behave} yield that there exists some
$\bar{y}$ and stratum $\mathcal{M}^{\prime}$ containing
$\bar{z}:=\left(  \bar{x},\bar{y}\right)  $ such that
$N_{\mathcal{M}^{\prime }}\left(  \bar{z}\right)
\cap\mathcal{L}^{\perp}\supsetneq\left\{ \mathbf{0}_{n+m}\right\} $.
Finally, since there are only finitely many strata,
Lemma~\ref{lemma:Ioffes-lemma} tells us that $S(x)$ and $\bar{S}(x)$
may differ only on a set of dimension at most $n-1$. This proves the
result in this particular case.\smallskip\newline We now consider
the case where $\mathcal{X}\neq\mathbb{R}^{n}$. Let
$\mathcal{X}=\dot{\cup}\mathcal{X}_{j}$ be a stratification of
$\mathcal{X}$, and let $\mathcal{D}$ be the union of strata of full
dimension in $\mathcal{X}$. Each stratum in $\mathcal{D}$ is
semialgebraically homeomorphic to $\mathbb{R}^{d}$, where
$d:=\dim\mathcal{X}$ and let
$h_{j}:\mathbb{R}^{d}\rightarrow\mathcal{X}_{j}$ denote such a
homeomorphism. By considering the set-valued maps $S\circ h_{j}$ for
all $j$, we reduce the problem to the aforementioned case. Since the
set of strata (\textit{a~fortiori} the set of full-dimensional
strata) is finite, we deduce that $S(x)\neq\bar{S}(x)$ can only
occur in a set of dimension at most $d-1$. $\hfill\Box$

\bigskip

The following result is now an easy consequence of the above.

\begin{theorem}
[Main result]\label{theorem:main-theorem} A closed-valued
semialgebraic set-valued map
$S:\mathcal{X}\rightrightarrows\mathbb{R}^{m}$, where
$\mathcal{X}\subset\mathbb{R}^{n}$ is semialgebraic, is strictly
continuous outside a set of dimension at most $\left(
\dim\mathcal{X}-1\right)  $.
\end{theorem}

\noindent\textit{Proof.} By Theorem~\ref{theorem:generic-osc} the map $S$
differs from the outer semicontinuous map $\bar{S}$ on a set of dimension at
most $\left(  \dim\mathcal{X}-1\right)  $. Apply
Theorem~\ref{theorem:generic-strict-cty}. $\hfill\Box$

\bigskip

\noindent\textbf{Remark.} Our main result
(Theorem~\ref{theorem:main-theorem}) as well as all previous
preliminary results
(Lemmas~\ref{lemma:decomposition},~\ref{lemma:normals-behave},
Theorems~\ref{theorem:generic-strict-cty},~\ref{theorem:generic-osc})
can be restated for the case where $S$ is definable in an o-minimal
structure. With slightly more effort we can further extend these
results in case where $S$ is tame, noting that one performs a
locally finite stratification in the tame case as opposed to a
finite stratification.

\section{Applications in tame variational analysis}

\label{sec-appl}

A standard application of Theorem~\ref{theorem:[RW98,5.55]} is to
take first the closure of the graph of $S$, and then deduce generic
continuity for the obtained set-valued map. While this operation is
convenient, this new set-valued map no longer reflects the same
local properties. For example, for a set $C\subset\mathbb{R}^{n}$,
consider the Hadamard normal cone mapping $\hat{N}_{C}:\partial
C\rightrightarrows\mathbb{R}^{n}$ and the limiting normal cone
mapping $N_{C}:\partial C\rightrightarrows\mathbb{R}^{n}$, where
$\mbox{\rm cl}(\mbox{\rm Graph}(\hat{N}_{C}))=\mbox{\rm
Graph}(N_{C})$. The Hadamard normal cone $\hat{N}_{C}\left(
\bar{z}\right)  $ for $\bar{z} \in\partial C$ depends on how $C$
behaves at $\bar{z}$, whereas the normal cone $N_{C}\left(
\bar{z}\right)  $ offers instead an aggregate information from
points around $\bar{z}$. The following result is comparable with
\cite[Proposition 6.49]{RW98}, and is a straightforward application
of Theorem~\ref{theorem:main-theorem}.

\begin{corollary}
[Generic regularity]Given closed semi-algebraic sets $C$ and $D$
with $D\subset C$, the set-valued map
$\hat{N}_{C}:C\rightrightarrows\mathbb{R} ^{n}$ is continuous on
$D\backslash D^{\prime}$, where $D^{\prime}$ is semialgebraic and
$\dim\left(  D^{\prime}\right)  <\dim\left(  D\right)  $. When
$D=\partial C$, we deduce that $\hat{N}_{C}\left(  z\right)
=N_{C}\left(  z\right)  $ for all $z$ in $\left(  \partial C\right)
\backslash C^{\prime}$, with $\dim\left(  C^{\prime}\right)
<\dim\left(
\partial C\right)$.
\end{corollary}

An analogous statement of the above corollary can be made for
(nonsmooth) tangent cones $\hat{T}_C$ and $T_C$ as well. \medskip

\noindent\textbf{Remark.} From the definition of subdifferential of
a lower semicontinuous function \cite[Definition~8.3]{RW98}, we can
deduce that the regular (Fréchet) subdifferentials are continuous
outside a set of smaller dimension. This result is comparable with
\cite[Exercise~8.54]{RW98}. Therefore nonsmoothness in tame
functions and sets is structured.

\bigskip

Let us finally make another connection to functions whose graph is a
finite union of polyhedra, hereafter referred to as \emph{piecewise
polyhedral functions}. Robinson \cite{Robinson81} proved that a
piecewise polyhedral function is calm (outer-Lipschitz) everywhere
\cite[Example~9.57]{RW98}, and a uniform Lipschitz constant suffices
over the whole domain of the function (although this latter is not
explicitly stated therein). A straightforward application of
Theorem~\ref{theorem:main-theorem} yields that piecewise polyhedral
functions are set-valued continuous outside a set of small
dimension. We now show that a uniform Lipschitz constant for strict
continuity applies.

\begin{proposition}
[Uniformity of graphical modulus] Let
$S:X\rightrightarrows\mathbb{R}^{m}$ be a piecewise polyhedral
set-valued map, where $X\subset\mathbb{R}^{n}$. Then $S$ is strictly
continuous outside a set $X^{\prime}$, with $\dim\left(
X^{\prime}\right)  <\dim\left(  X\right)  $. Moreover, there exists
some $\bar{\kappa}>0$ such that if $S$ is strictly continuous at
$\bar{x}$, then the graphical modulus $\mbox{\rm lip}_{X}S\left(
\bar{x}\mid\bar{y}\right)$ is a nonnegative real number smaller than
$\bar{\kappa}$.
\end{proposition}

\noindent\textit{Proof.} The first part of the proposition of strict
continuity is a direct consequence of
Theorem~\ref{theorem:main-theorem} since $S$ is clearly
semialgebraic. We proceed to prove the statement on the graphical
modulus. We first consider the case where the graph of $S$ is a
convex polyhedron. The graph of $S$ can be written as a finitely
constrained set $\mbox{\rm Graph}\left(  S\right)  =\left\{
z\in\mathbb{R}^{n+m}\mid Az=b,Cz\leq d\right\}  $ for some matrices
$A,C$ with finitely many rows. The projection of $\mbox{\rm
Graph}\left(  S\right)  $ onto $\mathbb{R}^{n}$ is the domain of
$S$, which we can again write as $\mbox{\rm dom}\left( S\right)
=X=\left\{  x\in\mathbb{R}^{n}\mid A^{\prime}x=b^{\prime},C^{\prime
}x\leq d^{\prime}\right\}  $. Let $\mathcal{L}$ be the lineality
space of $\mbox{\rm dom}\left(  S\right)  $, which is the set of
vectors orthogonal to the rows of $A^{\prime}$. We seek to find a
constant $\bar{\kappa}>0$ such that if $x$ lies in the relative
interior (in the sense of convex analysis) of $X$ and $y\in S\left(
x\right) $, then $\mbox{\rm lip}_{X}S\left( x\mid y\right)
\leq\bar{\kappa}$. By Proposition~\ref{pro:Mordukhovich-extended},
we have
\[
\bar{\kappa}=\sup_{\left(  x,y\right)  \in\mbox{r-int}\left(  X\right)
}\left\{  \frac{\left\vert a\right\vert }{\left\vert b\right\vert }\mid\left(
a,b\right)  \in N_{{\scriptsize \mbox{\rm Graph}\left(  S\right)  }}\left(
x,y\right)  \cap\left(  \mathcal{L}\times\mathbb{R}^{m}\right)  \right\}  ,
\]
where \textquotedblleft$\mbox{r-int}$\textquotedblright\ stands for
the relative interior. The above value is finite because of two
reasons. Firstly, if $\left( a,\mathbf{0}\right)  \in
N_{{\scriptsize \mbox{\rm Graph}\left( S\right)}}\left( x,y\right)
\cap\left( \mathcal{L}\times\mathbb{R}^{m}\right)  $ , then by the
convexity of $\mbox{\rm Graph}\left(  S\right)  $, $x$ lies on the
relative boundary of $X$. Secondly, the {}\textquotedblleft
sup\textquotedblright\ in the formula is attained and can be
replaced by {}\textquotedblleft$\max$\textquotedblright . This is
because the normal cones of $\mbox{\rm Graph}\left(  S\right)  $ at
$z=\left( x,y\right)  $ can be deduced from the rows of $C$ in which
$Cz\leq d$ is actually an equation, of which there are only finitely
many possibilities. In the case where $S$ is a union of finitely
many polyhedra, we consider the set-valued maps denoted by each of
these polyhedra. The maximum of the Lipschitz constants for strict
continuity on each polyhedral domain gives us the required
$\bar{\kappa}$. $\hfill\Box$

\vspace{1.0cm}

\noindent\textbf{Acknowledgement. } The majority of this work was performed
during a research visit of the second author at the CRM (Centre de Recerca
Matemàtica) in Barcelona (September to December 2008). The second author
wishes to thank his hosts for their hospitality.


\vspace{0.8cm}

\noindent Aris Daniilidis

\smallskip

\noindent Departament de Matemàtiques \newline Universitat Autònoma de
Barcelona \newline E-08193 Bellaterra, Spain

\smallskip

\noindent E-mail: \texttt{arisd@mat.uab.es} \newline\noindent
\texttt{http://mat.uab.es/\symbol{126}arisd}

\smallskip\noindent Research supported by the MEC Grant
MTM2008-06695-C03-03/MTM (Spain).

\bigskip

\noindent C. H. Jeffrey Pang

\smallskip

\noindent Center for Applied Mathematics \newline\noindent657~Rhodes Hall,
Cornell University, Ithaca, NY 14853, USA.

\smallskip

\noindent E-mail: \texttt{cp229@cornell.edu }\newline\noindent
\texttt{http://www.cam.cornell.edu/\symbol{126}pangchj}

\end{document}